# A Spiral Bicycle Track that Can Be Traced by a Unicycle

**Stan Wagon**, Macalester College, St. Paul, Minnesota, USA ⟨wagon@macalester.edu⟩

**Abstract.** A *unibike* curve is a track that can be made by either a bicycle or a unicycle. More precisely, the end of a unit tangent vector at any point on a unibike curve lies on the curve (so the bike's front wheel always lies on the track made by the rear wheel). David Finn found such a curve in 2002, but it loops around itself in an extremely complicated way with many twists and self-intersections. Starting with the polar square root curve $r = \sqrt{t/(2\pi)}$ and iterating a simple construction involving a differential equation apparently leads in the limit to a unibike curve having a spiral shape. The iteration gets each curve as a rear track of its predecessor. Solving hundreds of differential equations numerically, where each depends on the preceding one, leads to error buildup, but with some care one can get a curve having unibike error less than $10^{-7}$. The evidence is strong for the conjecture that the limit of the iteration exists and is a unibike curve.

## 1. The First Unibike

The two wheels of a rolling bicycle generally make two different tracks (Fig. 1). Because the vertical plane of the rear wheel always contains the point where the front wheel touches the ground (ignore the angled front fork of a real bicycle, which makes this only approximately true), the tangent to the rear wheel path at time $t$ strikes the front wheel path at the point where the front wheel is at time $t$. The distance between the two wheels for all bicycles here is one unit. Looking at several tangents generally allows one to determine which track is the front wheel as well as the direction of travel. See [**3**, **9**] for a discussion of tracks where the direction cannot be determined.

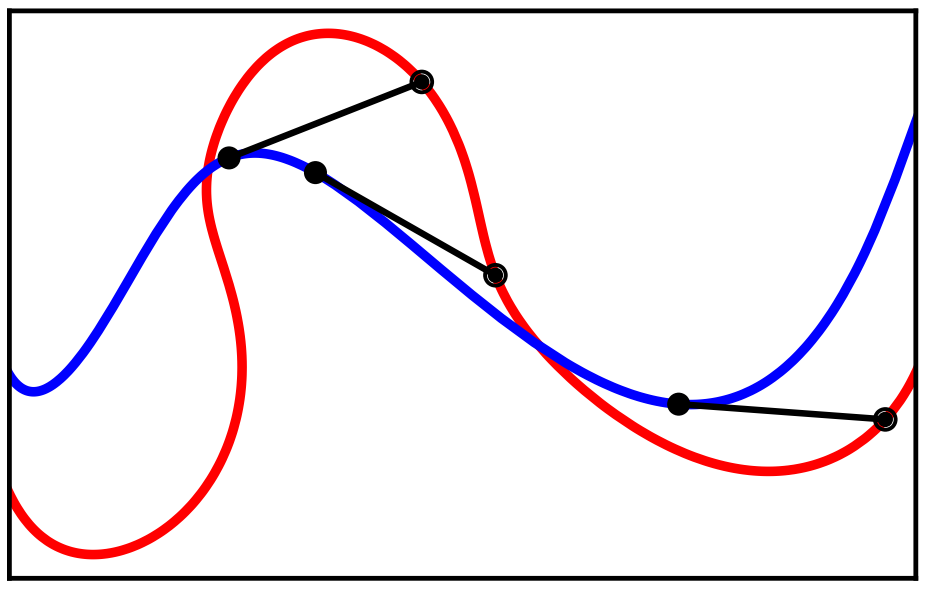

**Figure 1.** A set of front and rear tracks made by a bicycle. The unit tangent to the rear curve always ends on the front curve; this allows one to identify the rear track (blue) and the direction of travel (left to right).

A natural question is whether a bicycle can roll in such a way that the two wheels follow the same path. Of course, this happens if the path is a straight line, a case that we exclude. Here a curve is called *smooth* if it is continuously differentiable and the derivative is never $(0, 0)$.

**Definition.** Let $f(t)$ parametrize a smooth curve $C$. The *front track*, $\Phi(f)(t)$, is $f(t) + \frac{f'(t)}{\|f'(t)\|}$. The *unibike error* of $f$ at $t$, denoted $E(f)(t)$, is $\|\Phi(f)(t) - P\|$, where $P$ is the closest point to $\Phi(f)(t)$ lying on $C$. A *unibike curve* is a curve that is not a straight line and has zero unibike error at all points; i.e., each point $\Phi(f)(t)$ lies on $C$.

The term *unibike* refers to the fact that the curve can be traced by either a unicycle or both wheels of a bicycle. While the motion might go from $-\infty$ to a finite time, we generally care only about paths that start at a point and then have infinite length from that point, or perhaps are infinite in both directions. In 2002 David Finn ([**2**]; see also [**5**]) used a nonanalytic function to construct a unibike path (Fig. 2); the red curve over $[0, 1]$ that starts the construction is $f(t) = 4\, e^{-1/(t(1-t))}$, which is infinitely flat at its two ends; the 4 is inessential but clarifies the curve. Starting with $f$, form $\Phi(f)$, $\Phi(\Phi(f))$, and so on; the first four iterations are shown at right in Figure 2. In [**7**, §1.2], we showed how Finn's curve can be extended infinitely in the reverse direction as well, giving a path $g(t)$ that is defined on $(-\infty, \infty)$ and is such that $\Phi(g) = g$ (Fig. 3). For more on Finn's construction, including an animation, see [**8**].

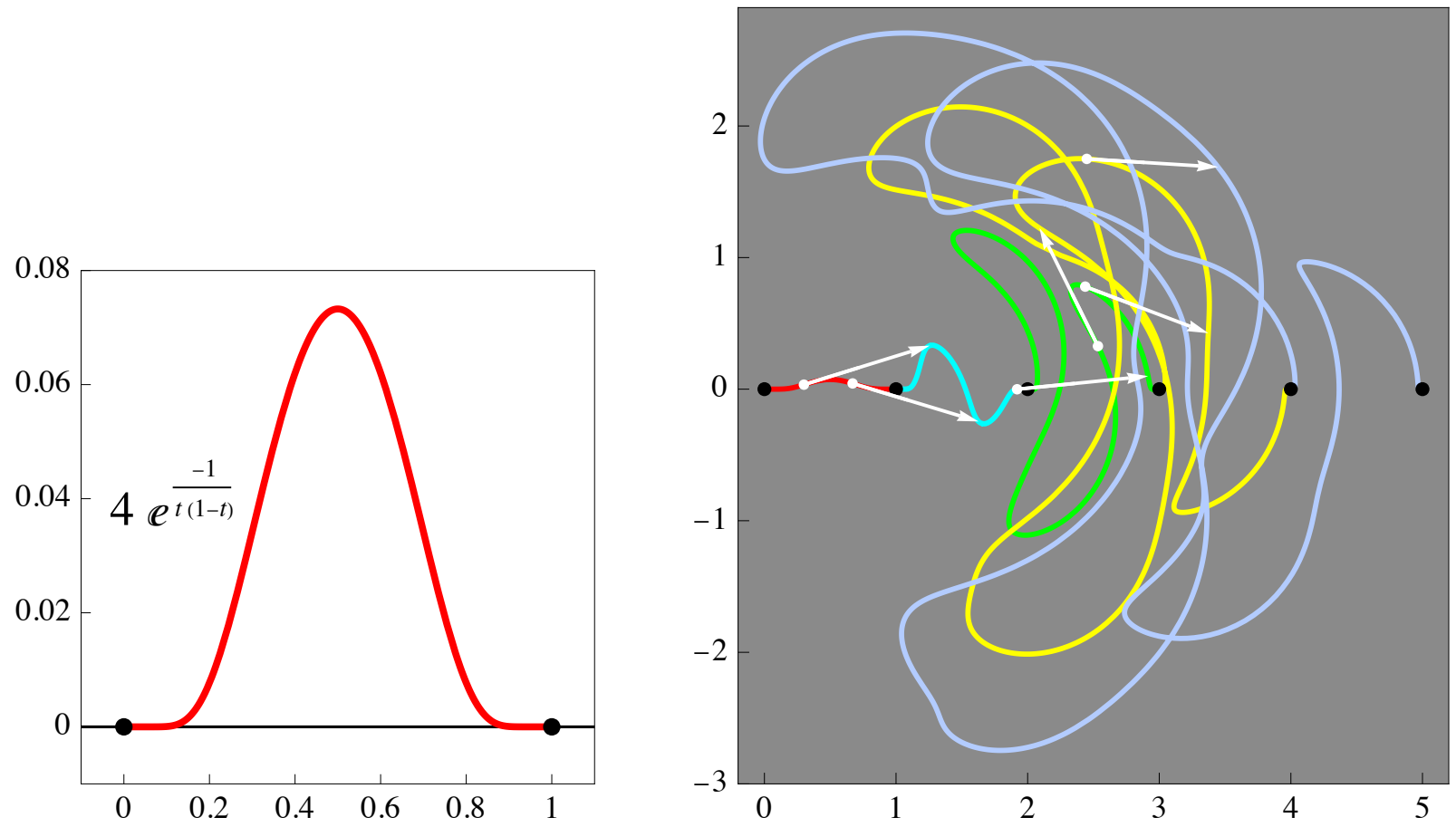

**Figure 2.** Finn's unibike path [**2**]. He started with the nonanalytic function at left. The repeated forward track of this curve gives the unibike curve at right, defined on $[0, \infty)$. The arrows show six bicycle positions. At each $(n, 0)$, $n \in \mathbb{N}$, the curve, viewed as a function of $x$, has infinitely many derivatives and they are all 0.

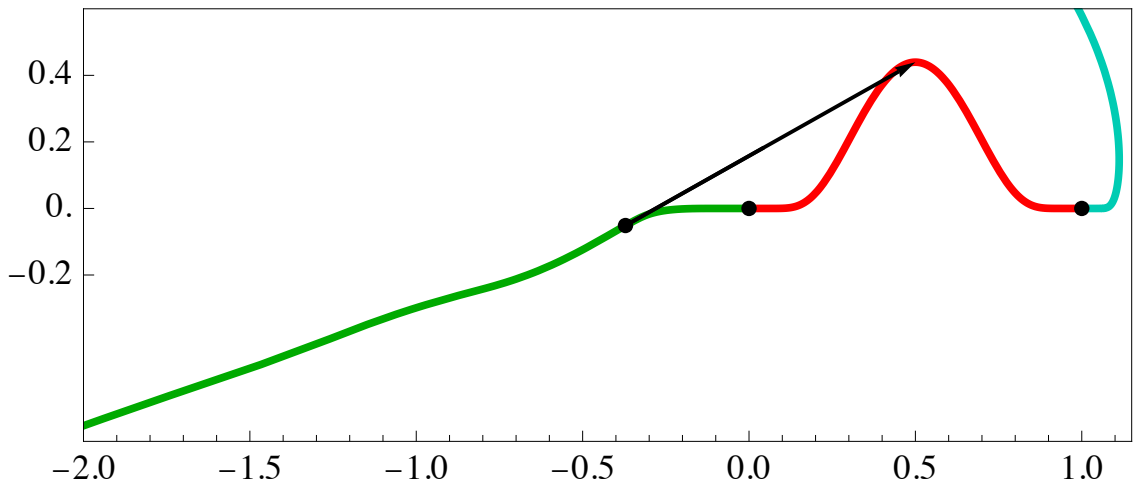

**Figure 3.** The green curve is part of a backward extension of Finn's track, obtained using the differential equation in §4. Here the seed curve (red) is $24\, e^{-1/(t(1-t))}$. The fact that it can be extended to $(-\infty, 0)$ is proved in [**7**], though it was not proved that the backwards curve has no cusps. Numerical work indicates that it is asymptotic to a straight line.

The goal of this paper is to present evidence for a spiral unibike curve that is infinite in one direction; it has no self-intersections. Figure 4 shows an example from §6 in which the front points are all within $10^{-6}$ of the rear track. Notably, this curve was obtained by solving only 11 differential equations numerically.

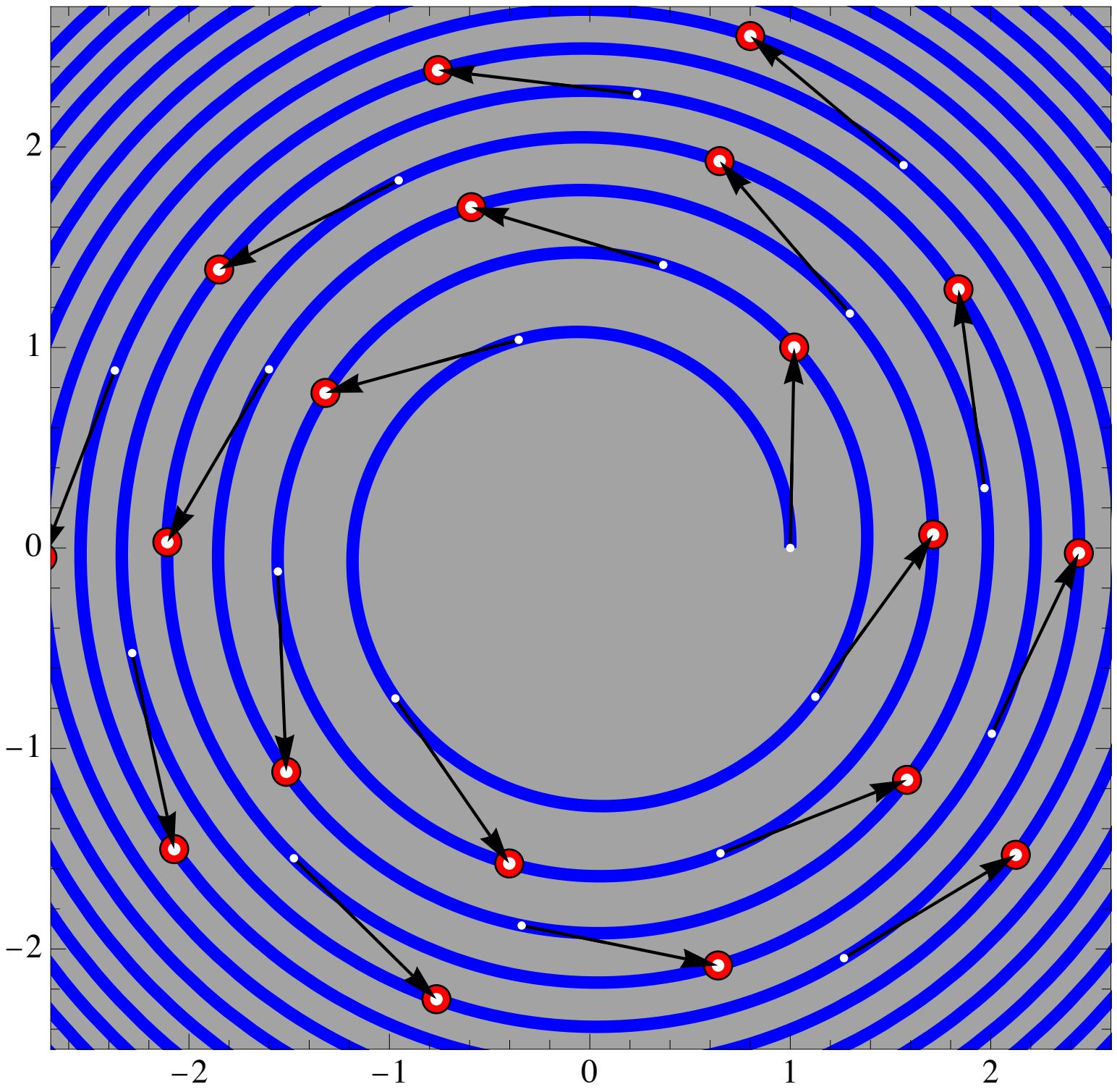

**Figure 4.** A nearly perfect unibike curve that is a weighted average of $K_6$ and $K_{12}$, defined in §6. The front track, defined by the ends of the unit tangent vectors, is never more than $10^{-6}$ units from the rear track (blue).

Just as Finn's example starts with a seed function—$e^{\frac{1}{t(1-t)}}$—a spiral unibike is the front-track extension of a seed curve. For the unibike shown in Figure 4, the seed curve is as in Figure 5.

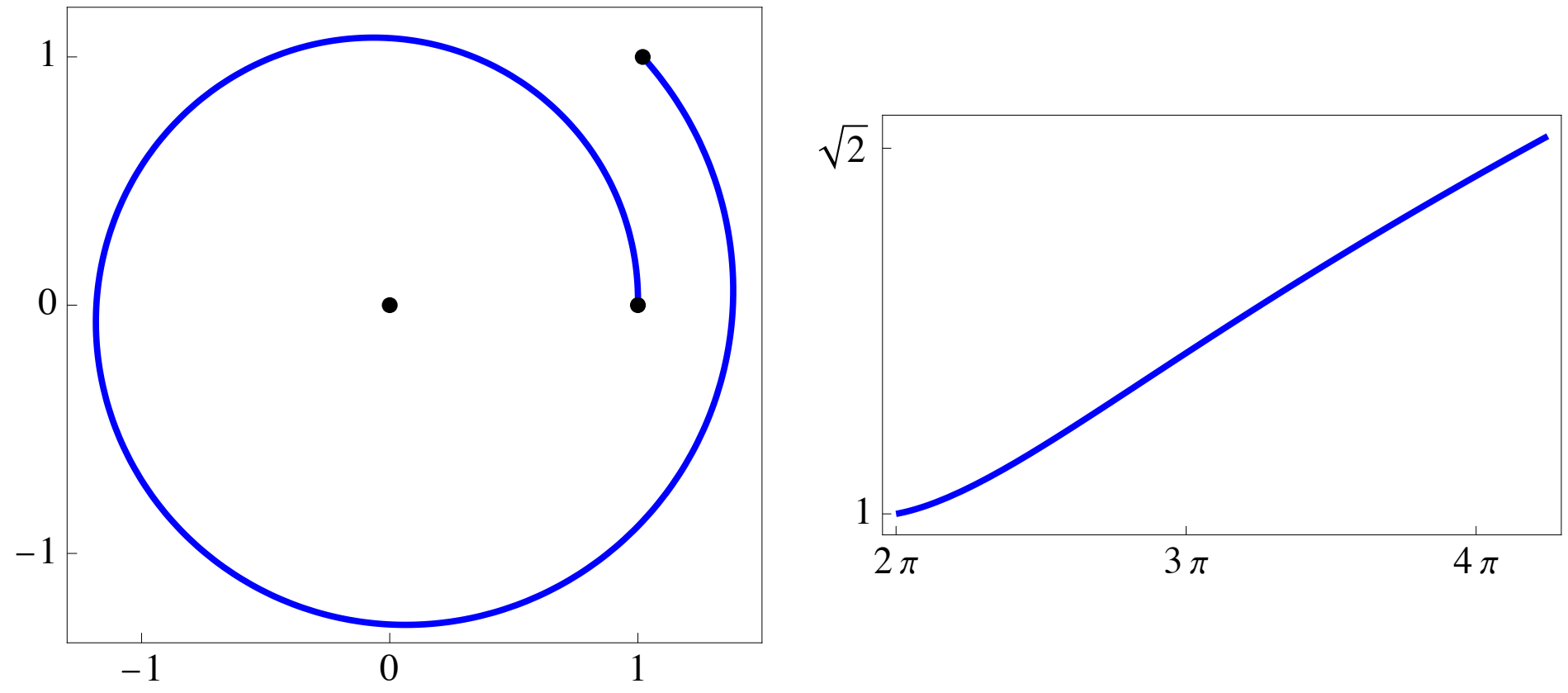

**Figure 5.** An approximation to a seed curve for a unibike. The polar radius is shown at right. That graph is sufficient to generate the full half-infinite unibike spiral.

The work here involves complicated symbolic expressions and *Mathematica* is used for all simplifications, limits, and series expansions, as well as for the intense numerical work. An Appendix contains the code for some of the key steps.

## 2. A Spiral Pseudo-Unibike

A natural idea is to search for a unibike curve that is a spiral (a curve with increasing polar radius, and hence having no self-intersections). In [**7**] we observed that the polar square root $r = \sqrt{\theta}$ comes close. Here we use the polar form $F_1(t) = \sqrt{\frac{t}{2\pi}}\,(\cos t, \sin t)$, starting at $(1, 0)$. This function is not analytic as $t$ rises to $\infty$. Figure 6 shows the curve $F_1$ as the rear wheel (blue) and the corresponding front wheel (yellow); we will use $F_0(t)$ for $\Phi(F_1)\,(t)$, but the parameter in $F_0$ is not the polar angle. Let $\alpha\,(t)$ be the proper angular parameter for $F_0(t)$; then $\alpha(t) = t + \beta(t) + 2\,\pi\,k$, where $0 \le \beta\,(t) < 2\,\pi$. The law of cosines on the triangle $((0, 0), F_1(t), F_0(t))$ yields

$$\beta(t) = \cos^{-1}\left(\frac{\sqrt{t}\,\gamma + \sqrt{2\,\pi}}{\sqrt{\gamma((t+2\,\pi)\,\gamma + 2\,\sqrt{2\,\pi\,t}\,)}}\right) = \tan^{-1}\left(\frac{2\,t}{1 + \sqrt{\frac{t}{2\pi}}\,\gamma}\right),$$

where $\gamma = \sqrt{1 + 4\,t^2}$ (Fig. 7); because $\gamma > 2\,t$, we have $0 < \beta\,(t) < \pi/4$. It is easily seen by comparing the norm of $F_0$ with that of $F_1(t + 2\,\pi + \beta(t))$ and $F_1(t + 4\,\pi + \beta(t))$ that the front point lies between the two mentioned points on $F_1$, which means that $k = 1$ and $\alpha(t) = t + 2\,\pi + \beta(t)$. An alternate form is $\alpha(t) = t + 2\,\pi + \mathrm{mod}\left[\tan^{-1}\left(F_0(t)_x, F_0(t)_y\right) - t,\ 2\,\pi\right]$, where the 2-argument form of arctan (*atn2* in FORTRAN; same as *arg* for complex numbers) is used. One can compute the inverse function $\alpha^{-1}(t)$ by numerical root-finding and then the parametric curve $F_0\left(\alpha^{-1}(t)\right)$ has the property that, for any $t$, the polar angle of $F_0\left(\alpha^{-1}(t)\right)$ is $t$.

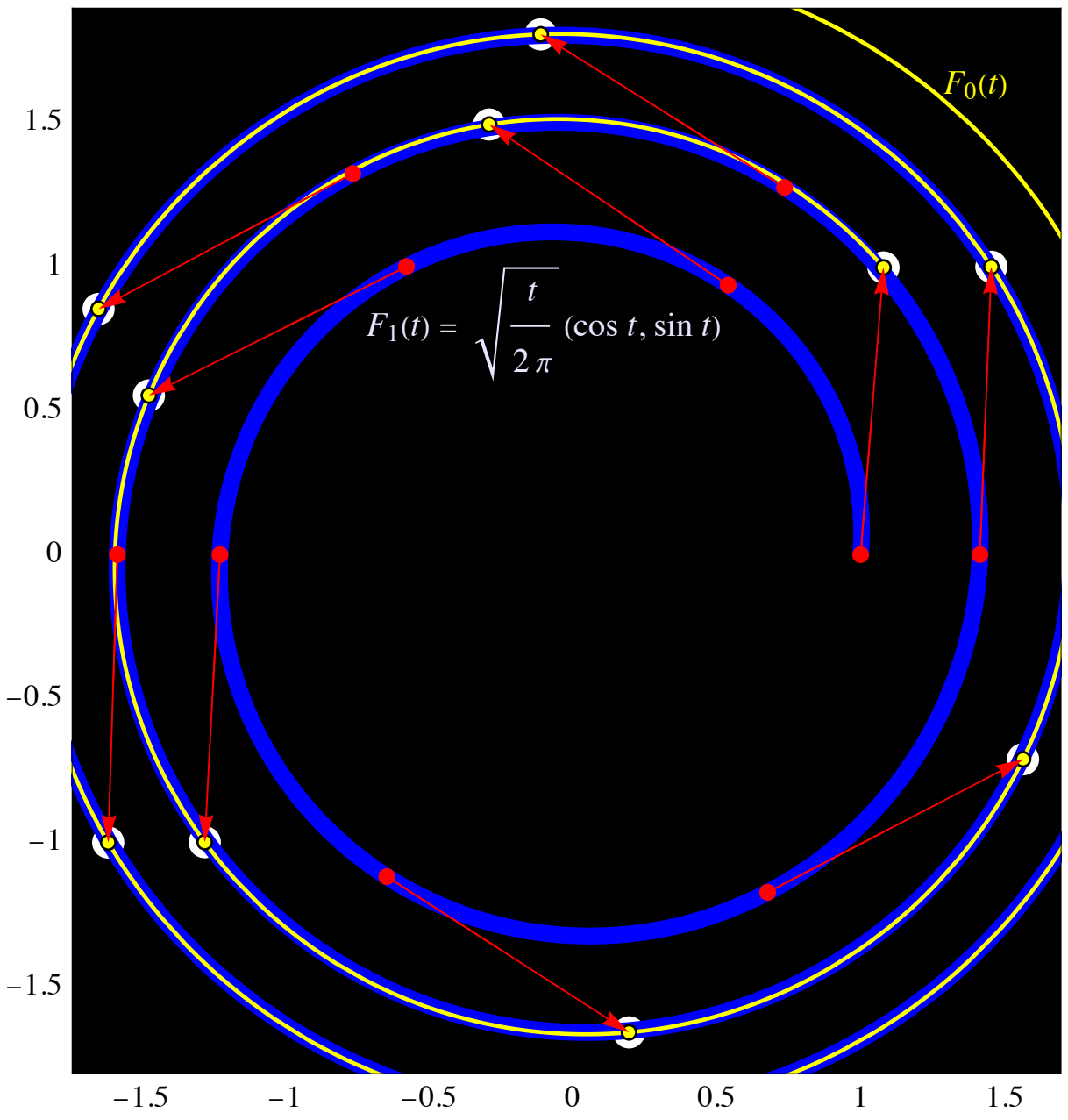

**Figure 6.** The polar square root spiral $F_1(t)$ (blue, with $2\,\pi \le t < \infty$) and its front wheel curve $F_0(t)$ (yellow). The small error (about $1/73$) is visible at the start—the end of the first tangent vector does not lie on the blue curve—but approaches zero quadratically as $t$ rises.

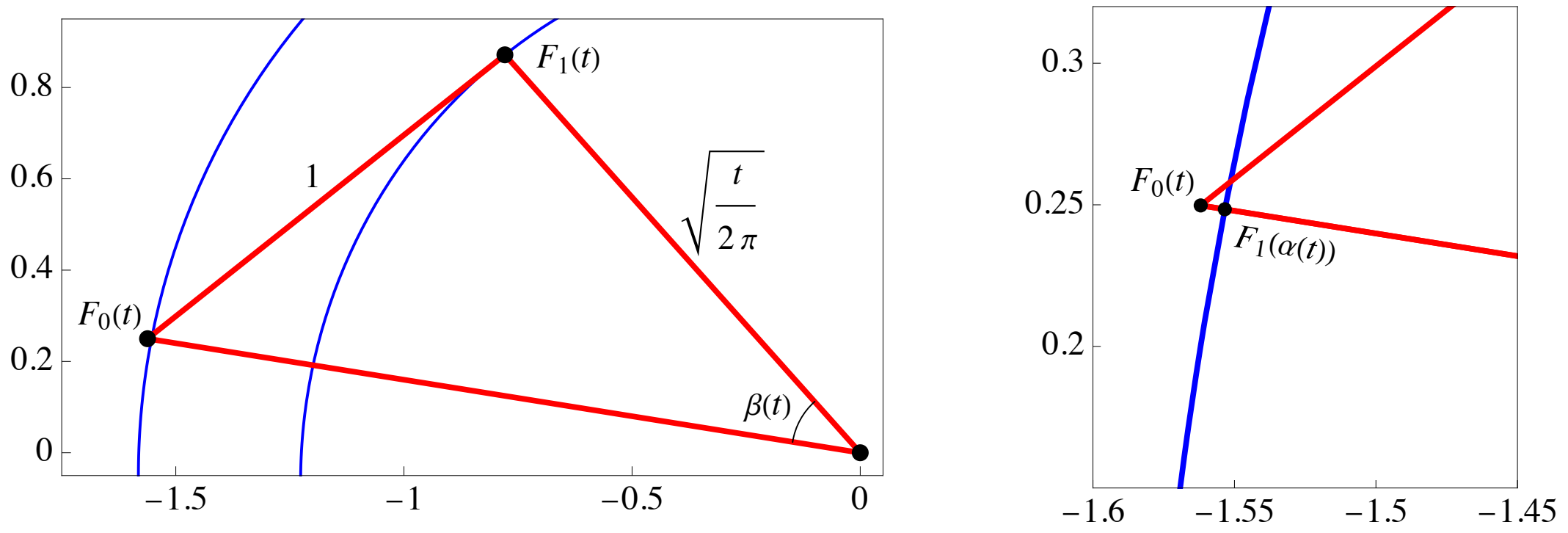

**Figure 7.** The triangle used to derive $\beta(t) = \tan^{-1}\left(\frac{2\,t}{1+\sqrt{t/(2\,\pi)}\;\gamma}\right)$. The image at right shows $F_1(\alpha(t))$, a good approximation to $F_0(t)$.

Our goal is to present strong evidence that a spiral unibike curve exists; more precisely, we seek a spiral unibike curve that starts at $(1, 0)$ and is close to the polar square root curve. Note that it is easy to get a spiral with arbitrary small error: just use the polar square root starting from a large value $t_0$; the error from that point on will be nicely small and as $t_0$ increases the curve's maximum unibike error approaches 0. But this process of deleting initial segments leads to no limiting curve, so is useless. Similarly, the polar curve $r = \theta^{1/n}$ for large $n$ has loops that are so close together that the maximum unibike error approaches 0 as $n \to \infty$; but again there is no limiting curve. The goal here is to construct a sequence of smooth curves $F_n(t)$ such that $F_n(2\,\pi) = (1, 0)$, $F_n$ has the same general shape as $F_1$, and $F_n$ converges to a spiral with zero unibike error.

We use $f \sim g$ to mean that $f$ is asymptotic to $g$ (the limit of the quotient is 1); when we have unproved but computationally evident closeness relations, we use $\approx$. Let $E_1(t) = E(F_1)\,(t)$ and more generally use $E_n$ for $E(F_n)$, where $F_n$ is defined in §4. It is clear that $E_1(t) \to 0$ as $t \to \infty$ because the loops get closer together. This simplistic view leads easily to $E_1(t) \leq \sqrt{\pi/8}\;t^{-1/2}$, but the true unibike error is much less. The end of the arrow in Figure 6 starting at $F_1(2\,\pi)$ is just 0.014 outside the blue curve; the simple bound gives $1/4$. A little farther along, the unibike error is invisible to the eye. Theorem 1 shows that $E_1(t) \sim (\pi/3)\,t^{-2}$.

This entire work requires a robust and fast method of computing the unibike error for a rear path $f\,(t)$. The generally best approach is the following. Let $P$ be the nearest point on $f\,(t)$ to a front track point $Q$; then the tangent to $f(t)$ at $P$ is perpendicular to $PQ$ and $P$ can be found by numerically finding $s$ such that $(f(s) - Q) \cdot f'(s) = 0$. However, it can happen that the derivative is difficult to compute (e.g., for $H_1$ in §6), and it that case a numerical minimum-finding method can be used.

For the polar square root one can quickly get the unibike error even when $t$ is $10^{10}$ or larger. The graph of $E_1(t)$ is in Figure 8; the log–log plot indicates that $\log(E_1(t)) - \log\left(\frac{\pi}{3}\,t^{-2}\right) \to 0$, which indicates that $E_1(t) \sim \frac{\pi}{3}\,t^{-2}$. This is proved in §3; an improved error model is $\frac{\pi}{3}\,t^{-2} + \sqrt{\frac{\pi}{8}}\;t^{-5/2} - \frac{11}{15}\,\pi^2\,t^{-3}$.

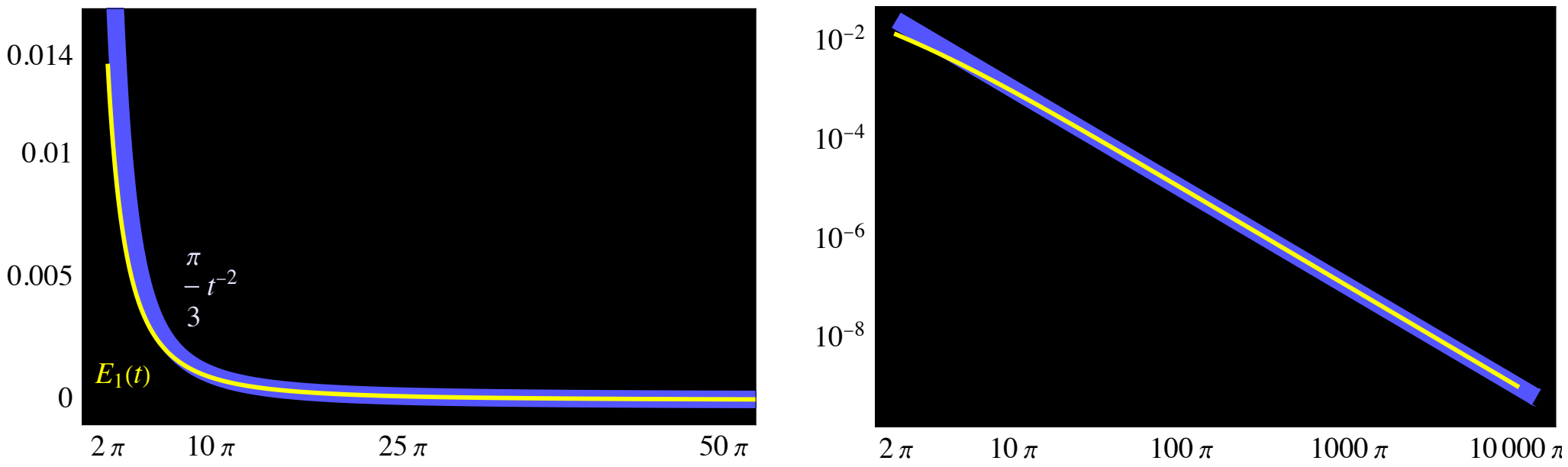

**Figure 8.** The polar square root error (left) and in log–log form (right). The error is asymptotic to $(\pi/3)\, t^{-2}$.

## 3. Error Bounds for the Polar Square Root

Let $C$ be the polar square root spiral defined by $F_1(t)$, $t \geq 2\pi$. Any point on $C$ provides an upper bound on the unibike error $E_1(t)$: the minimum distance from $F_0(t)$ to $C$. The point $F_1(\alpha(t))$, the first intersection of the line from $F_0(t_0)$ to the origin with $C$ (Fig. 7), is a good choice and yields a very tight asymptotic upper bound. The proof of a lower bound uses a geometric argument based on a tangent line to $C$.

**Notation.** Let $R_1(t) = \|F_1(t)\| = \sqrt{\frac{t}{2\pi}}$ and

$$R_0(t) = \|F_0(t)\| = \sqrt{1 + \frac{t}{2\pi} + \frac{1}{\gamma}\sqrt{\frac{2t}{\pi}}} = \sqrt{1 + \frac{t}{2\pi} + \frac{2}{\gamma} R_1(t)}\,.$$

Let $d((X, Y), Z)$ be the shortest distance from point $Z$ to the line through $X$ and $Y$. Let $\lambda(t) = \frac{\pi}{3} t^{-2} + \frac{1}{2\sqrt{2}} \pi^{1/2} t^{-5/2} - \frac{11}{15} \pi^2 t^{-3} - \frac{11}{6\sqrt{2}} \pi^{3/2} t^{-7/2}$.

**Theorem 1 (a).** $E_1(t) \leq R_0(t) - R_1(\alpha(t)) = \lambda(t) + \left(\frac{103}{70}\pi^2 - \frac{5}{8}\right)\pi\, t^{-4} + O\left(t^{-9/2}\right)$. The coefficient of $t^{-4}$ is about 43.66.

**(b)** $E_1(t)$ is asymptotically greater than or equal to $d[(F_1(\alpha(t)), F_1(\alpha(t)) + F_1'(\alpha(t))), F_0(t)]$, which equals $\lambda(t) + \left(\frac{103}{70}\pi^2 - \frac{2}{3}\right)\pi\, t^{-4} + O\left(t^{-9/2}\right)$. The coefficient of $t^{-4}$ is about 43.53.

**(c)** $E_1(t) = \frac{\pi}{3} t^{-2} + \frac{1}{2\sqrt{2}} \pi^{1/2} t^{-5/2} - \frac{11}{15}\pi^2 t^{-3} - \frac{11}{6\sqrt{2}} \pi^{3/2} t^{-7/2} + O\left(t^{-4}\right)$.

**Proof. (a)** For $t \geq 2\pi$, let $P = F_0(t)$ and let $\delta(s)$ be the distance from $P$ to $F_1(s)$. For any $s$, $E_1(t) \leq \delta(s)$. Let $Q = F_1(\alpha(t))$; $Q$ is on the line from $(0, 0)$ to $P$ (Fig. 6, right) and $\|Q\| = \sqrt{\frac{\alpha(t)}{2\pi}}$. We then have

$$\delta(\alpha(t)) = \|P\| - \|Q\| = R_0(t) - \sqrt{\frac{s}{2\pi}} = \frac{1}{\sqrt{2\pi}}\left(\sqrt{1 + \frac{t}{2\pi} + \frac{1}{\gamma}\sqrt{\frac{2t}{\pi}}} - \sqrt{\alpha(t)}\right),$$

and a series expansion around $t = \infty$ using the formula for $\alpha(t)$ gives the claimed result. Code for the series is in the Appendix.

**(b)** The equality is the 4th-order Taylor polynomial at $t = \infty$ applied to the standard point-to-line distance formula. For the inequality, let $P = F_0(t)$ and let $L$ be the tangent line to $\mathcal{S}$, the curve defined by $F_1$, at $Q = F_1(\alpha(t))$ (Fig. 9, left). Then for any point $W$ on the near side of $L$ (the side containing the origin), $\|P - W\| \geq d(L, P)$.

To finish, we show that the closest point on $C$ to $P$ cannot lie on the far side of $L$. It suffices to consider points $F_1(s)$, where $t + 4\pi \le s \le F_1(\alpha(t) + 2\pi)$ (red in Fig. 9). The lower bound is justified by the fact that $L$ cannot strike $C$ at $F_1(s)$ where $s \le t + 4\pi$. This is proved by the standard dot-product formula to show that $F_1(t + 4\pi)$ is right of the segment from $Q$ to $Q - F_1'(\alpha(t))$; the dot product is asymptotic to $-\pi / t$, and so is asymptotically negative, as desired (code for this symbolic work is in the Appendix). And the upper bound is justified because of the radial expansion of $C$: points beyond $F_1(t + 4\pi)$ can be shown to be too far from $P$ by reflection in the polar line through $P$ and using the main case. For the main case, we need only observe that the minimum distance in the domain is greater than the minimum distance from $P$ to the line segment joining $F(\alpha(t) + 2\pi)$ and $F(t + 4\pi)$ (the green line in Fig. 9). That distance, by the standard point-to-line distance formula, is asymptotic to $\sqrt{\pi/2}\, t^{-1/2}$, vastly greater than $(\pi/3)\, t^{-2}$, the bound from (a).

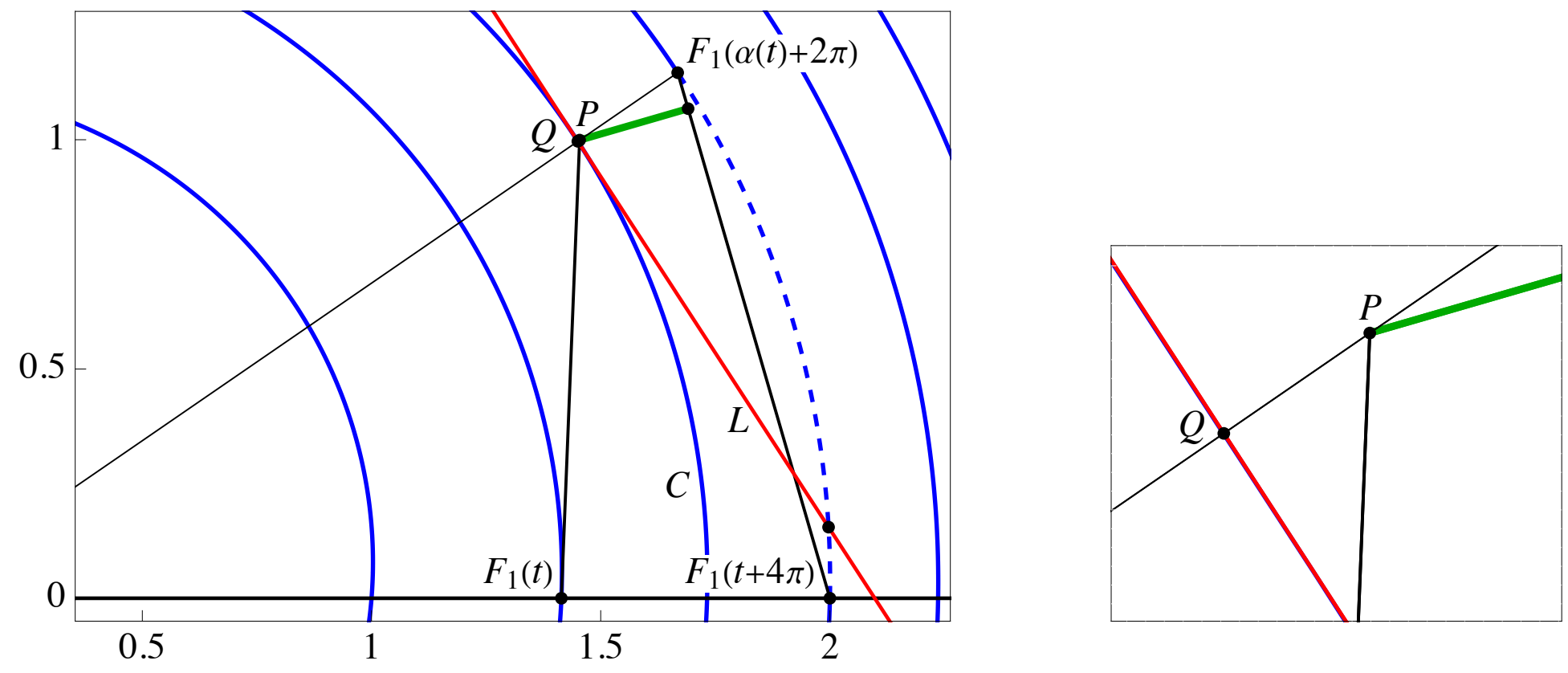

**Figure 9.** Here $P = F_0(t)$, $Q = F_1(\alpha(t))$, and $L$ is tangent to $C$ at $Q$. Any point on the origin's side of $L$ is farther from $P$ than $Q$ is. The tangent passes through $C$ only beyond $F_1(t + 4\pi)$. The shortest distance from any point on the dashed part of $C$ cannot be less than the shortest distance from $P$ to the black line. Here $t = 4\pi$. A close-up is at right.

**(c).** This follows immediately from (a) and (b). □

The approximations of Theorem 1 are very good. When $t = 10^6\, 2\pi$, a million loops of the curve, the two bounds differ by about $10^{-28}$ and the true value is within about $10^{-29}$ of the lower bound. Figure 10 shows how the lower bound differs from the true error by about $65\, t^{-9/2}$. The assertions of the theorem appear to be true absolutely, not just asymptotically, but the various expressions are too complicated to lead to immediate proofs of absolute results.

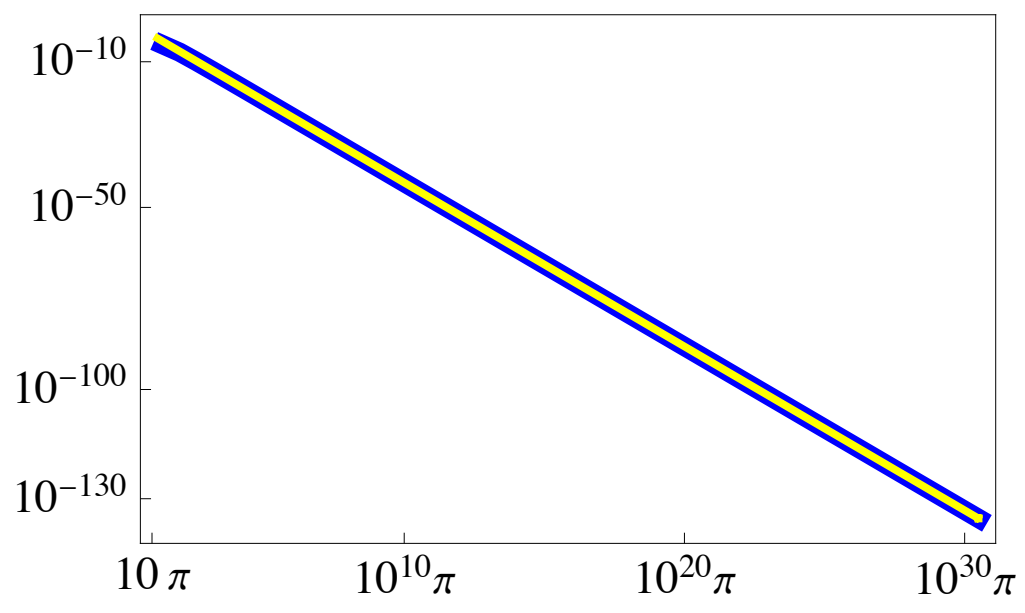

**Figure 10.** The difference between the polar square root and the lower bound of Theorem 1(b) is shown in blue. It is modeled closely by $65\, t^{-9/2}$, shown in yellow.

## 4. The Rear Track

Because the unibike error for the polar square root is so small, one has an irresistible urge to try to modify the curve so as to reduce the error. The hope is to find an iterative process to get curves $F_n$, each starting at $(1, 0)$, for which the unibike error approaches 0 and the curves approach a limiting curve, $F_\infty$. Then $F_\infty$ will be a unibike curve and, one hopes, in the shape of a spiral. There might well be no closed form for such a zero-error spiral, but the work here is aimed at showing that it is likely that a spiral unibike curve exists.

To get the front track $F_0$ from the rear track, just add the unit tangent vector. Figure 6 shows the front curve for the polar square root. Iterating this process is not useful because a loop is lost at each step and there is no limiting curve. But it is of interest to look at the unibike error for $F_0$: it is worse than $F_1$. For example, $E_1(\alpha(2\pi))$, the error for the polar square root near $F_0(2\pi)$, is 0.0043. The error at $F_0(2\pi)$ is roughly double at 0.0095. This makes the reverse direction seem promising: consider $F_1$ to be the front wheel track and look at the corresponding rear wheel track in the hope that this will significantly reduce the unibike error.

The path of the rear wheel corresponding to a given front wheel track can be obtained from a simple differential equation. The equation was first derived by Dunbar et al [**1**] and is an example of a Riccati equation (see [**6**, chap. 18], or [**7**]). Suppose $F(t)$ is the path of the front wheel and we seek the rear-wheel path $R(t)$. This task is the physically natural one because a bicyclist steers the front wheel and the rear wheel follows according to geometry. This idea applies also to cars, trucks, and buses, where again the steering is at the front and the rear follows. Because the rear wheel cannot steer, its velocity vector always points to the location of the front wheel. We can decompose the velocity vector of the front wheel into its components parallel to and perpendicular to the bike (Fig. 11) and only the parallel component affects the speed of the rear wheel. Therefore the speed of the rear wheel is the magnitude of the parallel component, which is $F'(t) \cdot (F(t) - R(t))$, yielding

(1) $$R'(t) = F'(t) \cdot (F(t) - R(t)) \, [F(t) - R(t)].$$

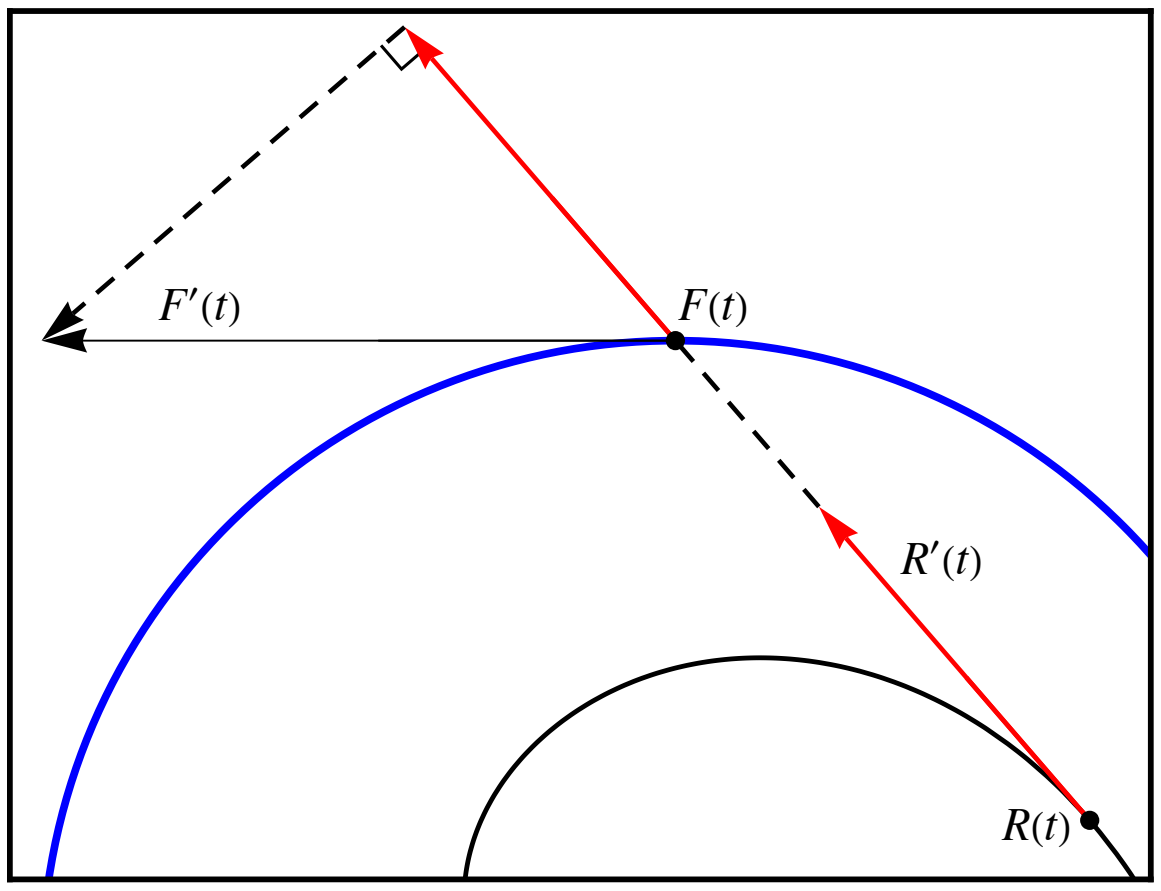

**Figure 11.** The instantaneous speed of the rear wheel is the magnitude of the component of $F'(t)$, the front wheel's velocity vector, in the direction the rear wheel must travel.

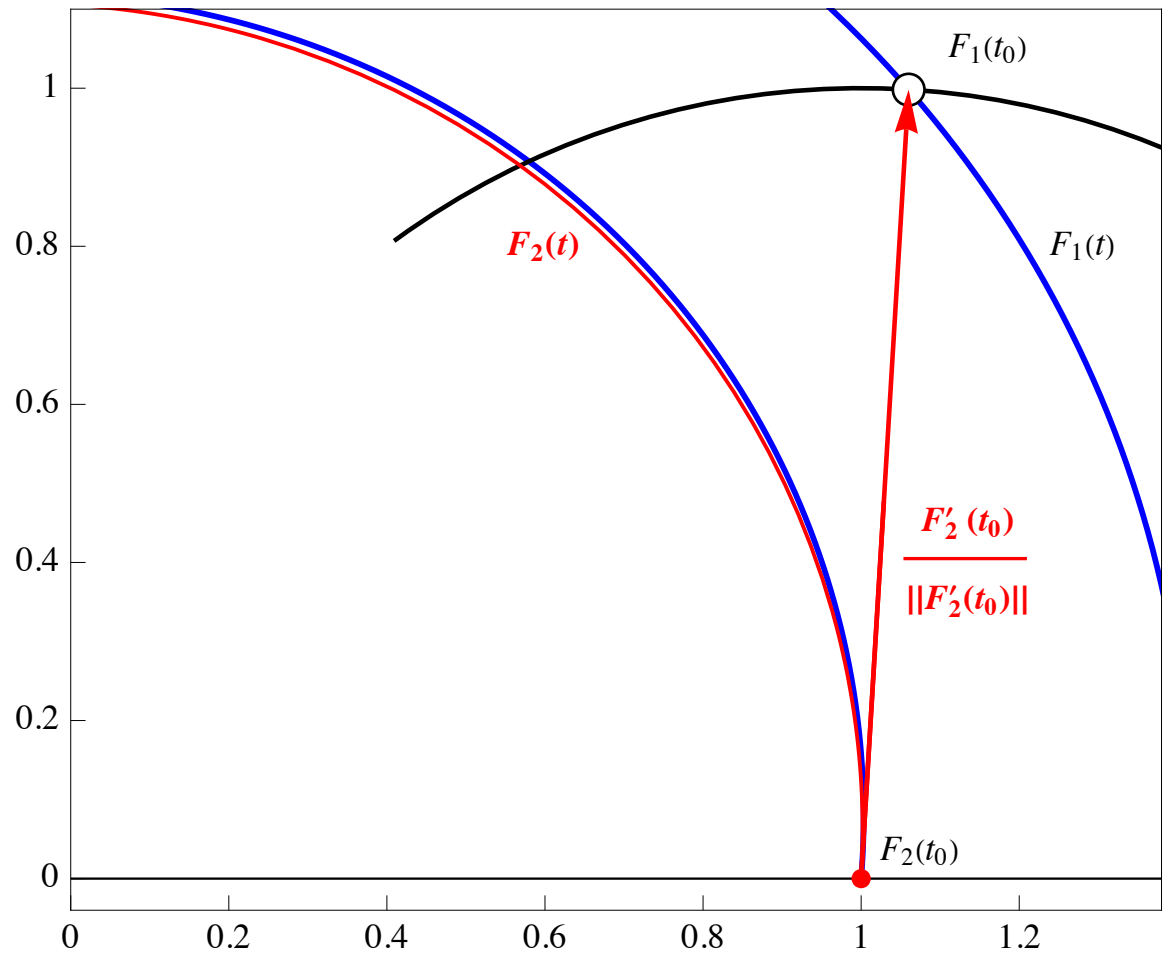

**Figure 12.** Using $t_0 = 13.3217\ldots$, the white disk, at $F_1(t_0)$, is at unit distance from $(1, 0)$. The initial condition for the differential equation to get the rear track $F_2$ from the polar square root $F_1$ is therefore $F_2(t_0) = (1, 0)$. The red curve is the differential equation solution giving $F_2$, and the red arrow confirms that the bicycle's front wheel is at $F_1(t_0)$ when the rear is at $(1, 0)$.

We wish to study the rear-wheel track, call it $F_2$, when the front track is the polar square root $F_1$. But this is not unique as it depends on the initial condition used. We want $F_2$ to start at $(1, 0)$, so we use the initial condition $F_2(t_0) = (1, 0)$ for an appropriate $t_0$; this means that $F_1(t_0)$ should be at distance exactly 1 from $(1, 0)$ and also, to preserve the spiral shape, $t_0$ should be the first choice greater than $4\pi$ (see Fig. 12). This value is easily found by root-finding: $t_0 \approx 13.3217$; the seed for this can be $4\pi$ or, better, $\alpha(2\pi)$, which is about 13.312. One cannot use the simpler condition $F_2(2\pi) = (1, 0)$ as the initial value because the equation $R'(t) = F'(t) \cdot (F(t) - R(t))\,[F(t) - R(t)]$ ties the front and rear together at the single time $t$. If $h(t)$ is the function such that the front wheel for $F_2(t)$ is at $F_1(h(t))$, then we have no information about $h(t)$ except (assuming $F_2(2\pi) = (1, 0)$) that $h(2\pi) = t_0$.

We can numerically solve the initial-value problem given by the rear-track differential equation and $F_2(t_0) = (1, 0)$ over $[t_0, t_{max}]$ to get $F_2$, a curve that starts at $(1, 0)$ and has $F_1$ as its front track. *Mathematica*'s numerical approach to differential equations gives the solution as an interpolating function over the requested domain, using piecewise cubics with matching derivatives at the junctions (Hermite interpolation). Complete code to get the rear track is only four lines (see Appendix). Figure 13 shows the two tracks along with some unit tangent vectors. Note that $F_2$ refers to the theoretical rear track defined for $2\pi \le t < \infty$, but we will use it also for a numerical approximation to the theoretical function. We use $E_2(t)$ for the unibike error of $F_2$ at $t$.

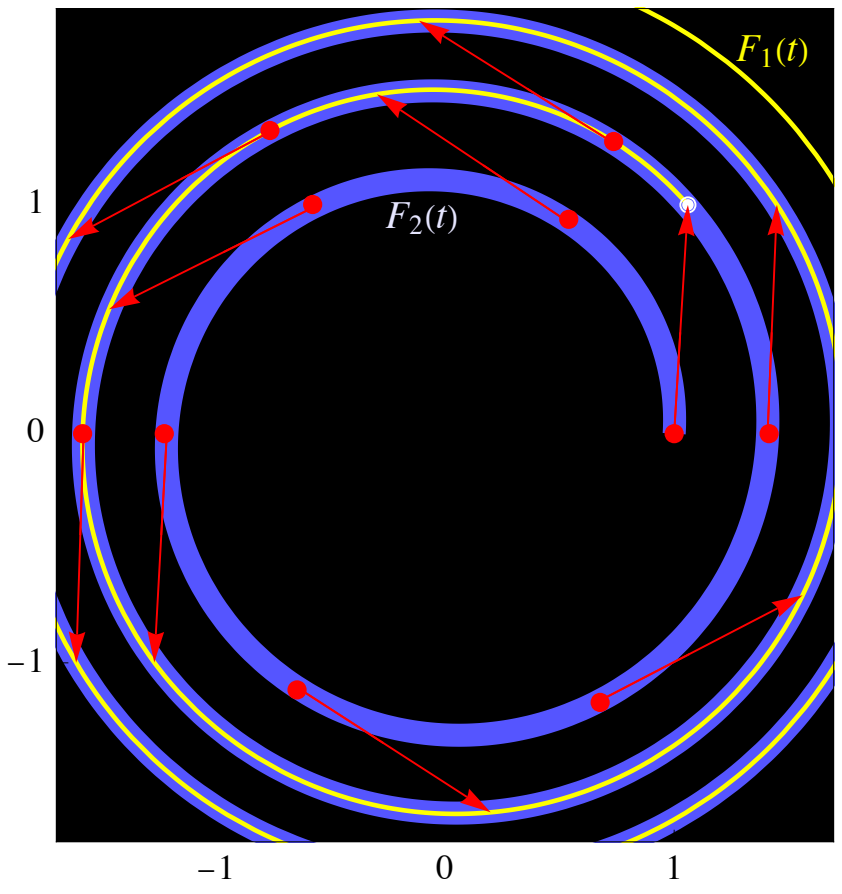

**Figure 13.** The polar square root $F_1$ is in yellow, with its rear track in blue and some arrows showing bicycle positions. Because $F_1$ is almost a unibike track, $F_1$ and $F_2$ almost coincide. The error $E_2(2\pi)$ is barely visible as the white circle misses the blue curve by less than 0.006.

It's always good practice to check the results of an approximate numerical solution. That can mean several things: (1) the difference between the left and right sides of the equation; (2) the difference between the computed solution and the true solution, or (3) the minimal distance between the front curve determined by the solution and the front curve used in the differential equation. We cannot get a symbolic solution, but we can compute a very high-precision solution and treat that as if it were the true solution. Figure 14 shows the errors in a test case up to $10\pi$; the errors are about what one would expect given the default working precision (16 digits) and precision goal (6 digits). The error here decreases when the working precision is increased, as is necessary when we iterate the process, but the computation time increases. A working precision of 22 and precision goal of 15 suffices for the computa tions in this project. For a domain out to $650\pi$, this leads to an interpolating function for $F_2$ that has more than 350,000 interpolating points.

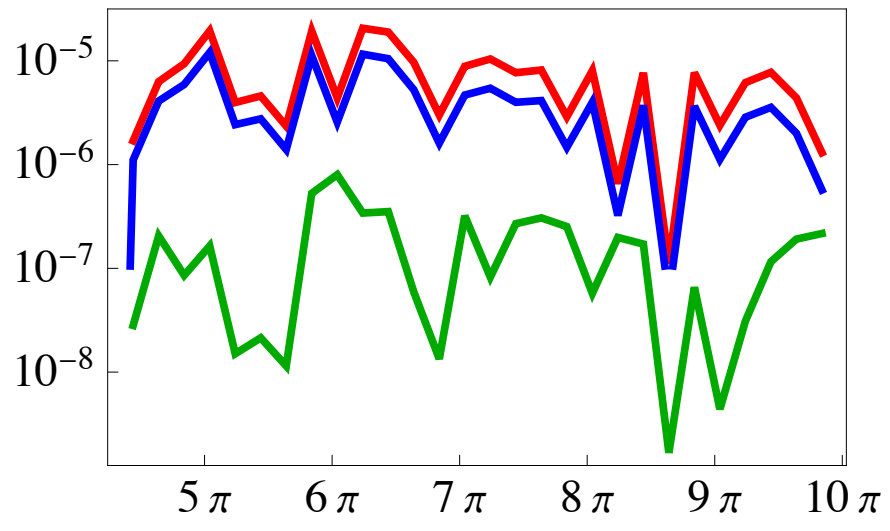

**Figure 14.** Here $F_2$ is the rear track obtained from $F_1$ using machine precision, 16 digits. The red curve indicates the difference between the two sides of the differential equation that defines $F_2$ from $F_1$. The blue points give the minimum distance from $\Phi(F_2)(t)$ to the curve $F_1$. The green values are the difference between the computed solution and a proxy for the true solution computed with 37 digits of working precision.

For several reasons (the error computation; iterating the process) it is important to reparametrize the solution to the differential equation so that the parameter is polar angle; that is, the parameter is the angle made by the point with $(0, 0)$ and the $x$-axis, increased by the correct integer multiple of $2\pi$. This is a little delicate to program in general because of several numerical issues that can arise. From now on we consider $F_2$ (and all $F_n$) to be so parametrized. Figure 15 compares the polar radii of $F_1$ and the reparametrized $F_2$.

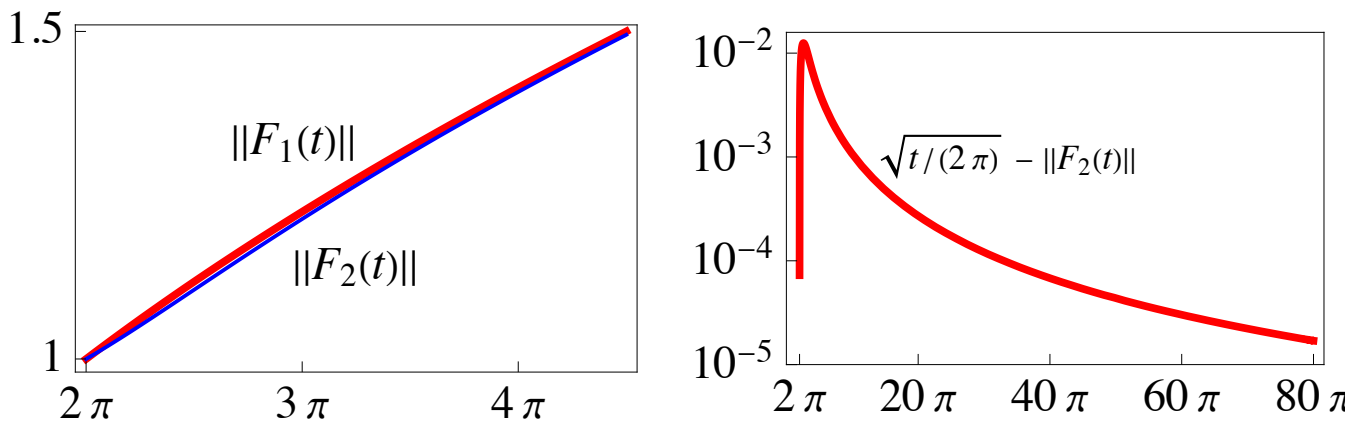

**Figure 15.** Left: The norms of $F_1$ and $F_2$ after reparametrization. Right: A log plot of the radial difference between $F_1$ and $F_2$. They both start at $(1, 0)$ which explains the spike near $2\pi$ as the difference can only rise from 0 before it starts to decrease back to 0.

Figure 16 is a log–log plot comparing the unibike error for $F_2$ with that for $F_1$. Theorem 1 proves that the asymptotic slope of the log–log form of $E_1(t)$ is $-2$. The asymptotic behavior is apparently unchanged when we move to $F_2$, but the overall error has decreased a lot. A key point is that the error $E_2$ is decreasing, so the largest unibike error occurs at $t = 2\pi$. Here we have $E_1(2\pi) = 0.0137\ldots$ and $E_2(2\pi) = 0.0058\ldots$; the maximum unibike error has been halved! In the next section we will see that the error reduction continues as the process is iterated, following a wonderfully simple pattern.

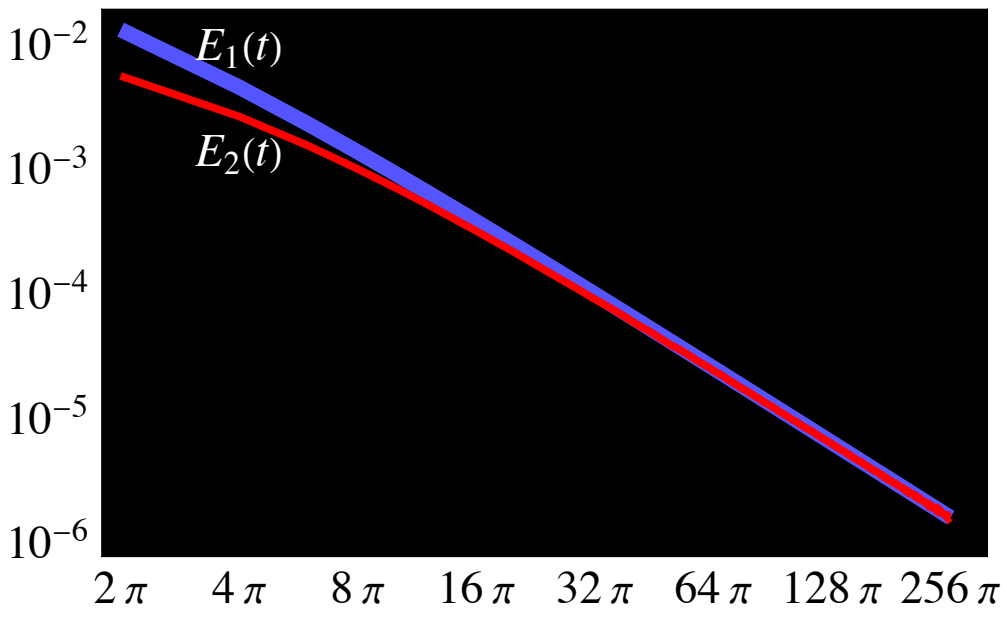

**Figure 16.** A logarithmic view of the unibike error for $F_1$ and $F_2$.

## 5. A Spiral Unibike

Because $E_2$ is less than $E_1$, one can hope that repeating the process will yield a sequence $F_n$ such that the corresponding unibike errors $E_n$ will continually decrease. This would not mean they go to zero, so the hope is to find a definitive pattern in the errors that implies a limit of 0. There is indeed such a pattern! In §6 we will present a heuristic argument to support the experimental result that $E_n(2\pi) \approx \left(9\pi^{3/2}\right)^{-1} n^{-3/2}$. Because computations indicate that $E_n(t)$ is decreasing in $t$, and further computa tions and another heuristic argument indicate that the sequence $F_n$ approaches a limiting function $F_\infty$, this means that there is almost certainly a spiral unibike curve. In short, two miracles occur when we study the rear tracks: the unibike error decreases in a patterned way, and the tracks converge to a limiting track, again following a simple pattern.

With $F_{256}$ as the target, we must start with $F_2$ on a very large domain ($F_1$ is defined for all $t \geq 0$), as a loop is lost whenever we get a new rear track. Starting with $F_2$ defined to $t = 650\pi$, we can iterate the rear track computation of §4 (using a working precision of 22 decimal digits) and reparametrize at each step, to get $F_n$ for $2 \leq n \leq 256$, and with $F_{256}$ defined up to $133\pi$. This takes seven hours. Figure 17 shows $F_{256}$, along with a closeup view of its unibike error at $2\pi$, which is under $5 \cdot 10^{-6}$. We use $E_n(t)$ for $E(F_n)(t)$ and $R_n(t)$ for $\|F_n(t)\|$.

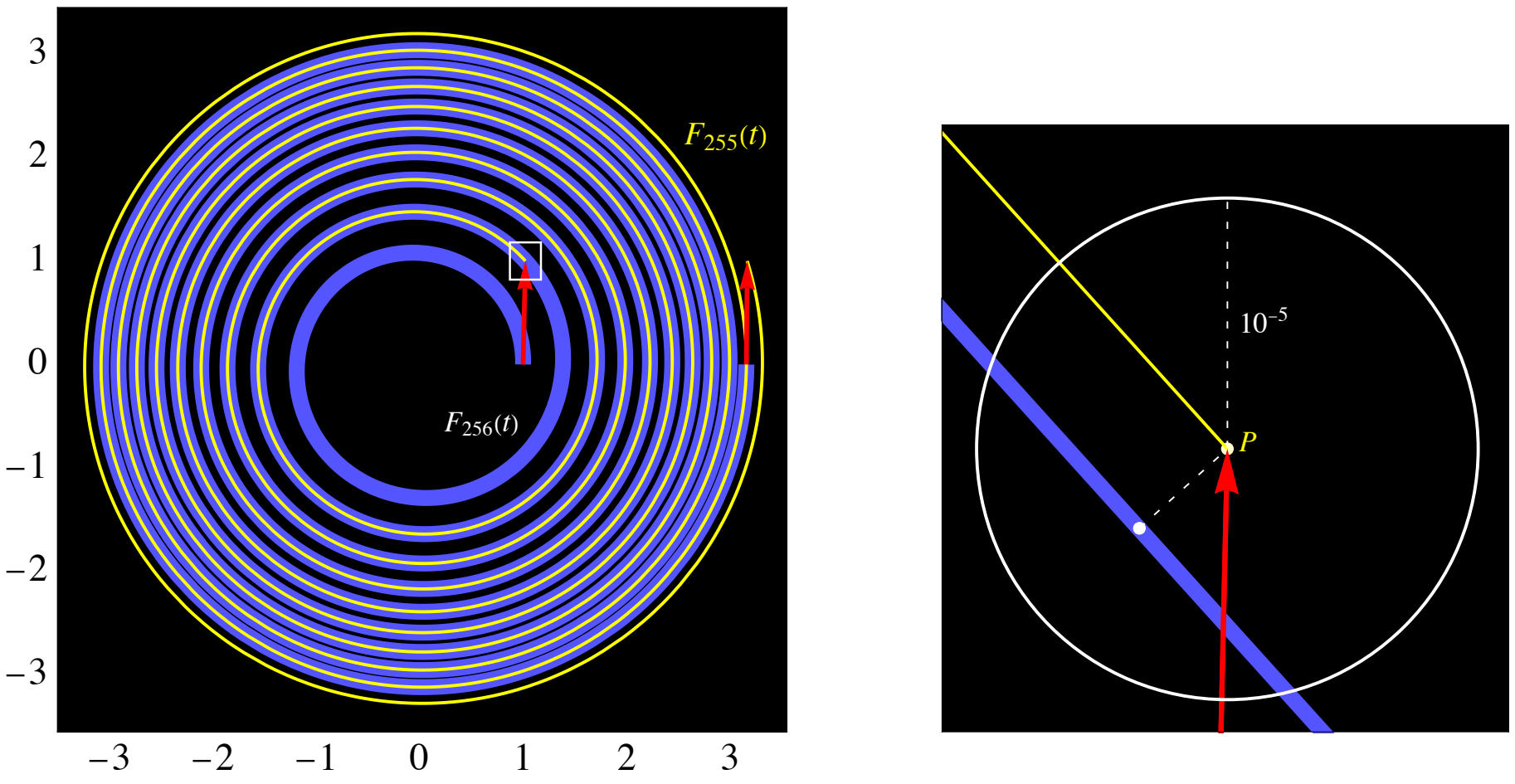

**Figure 17.** The track for $F_{256}$ (blue) and its front wheel $\Phi(F_{256})$, which equals $F_{255}$. At normal scale (left), $F_{256}$ appears to be a perfect spiral unibike: the maximum error is under $10^{-5}$. A zoom to the initial front point $P = \Phi(F_{256})(2\pi)$ is at right; the white point is where $F_{256}$ is closest to $P$.

A log–log plot of the error of $F_n$ at $2\pi$ (yellow dots in Fig. 18) shows the first miracle. The error decreases in a definite pattern as the essentially linear plot indicates a power rule. A fit gives the correct power as $n^{-3/2}$; for example, $E_{300}(2\pi) = 3.72\ldots \cdot 10^{-6}$ while $1 \big/ \left(9\,(\pi\, 256)^{3/2}\right) = 3.84\ldots \cdot 10^{-6}$. A subtle analysis in §6 gives evidence for $E_n(2\pi) \sim \frac{1}{9\,\pi^{3/2}}\, n^{-3/2}$. If, as expected, the error monotonicity for rising $t$ holds, this asymptotic result implies that $\lim_{n\to\infty} E_n(t) = 0$ for any $t$. Figure 19 shows the solid evidence for decreasing $E_n(t)$ as $t$ rises, and it appears that each $E_n(t) \sim E_1(t)$. The error model indicates that $E_{736} < 10^{-6}$, but computing such a distant $F_n$ is not possible. We will see in §6 how to get such small error using only 11 steps of the rear track construction.

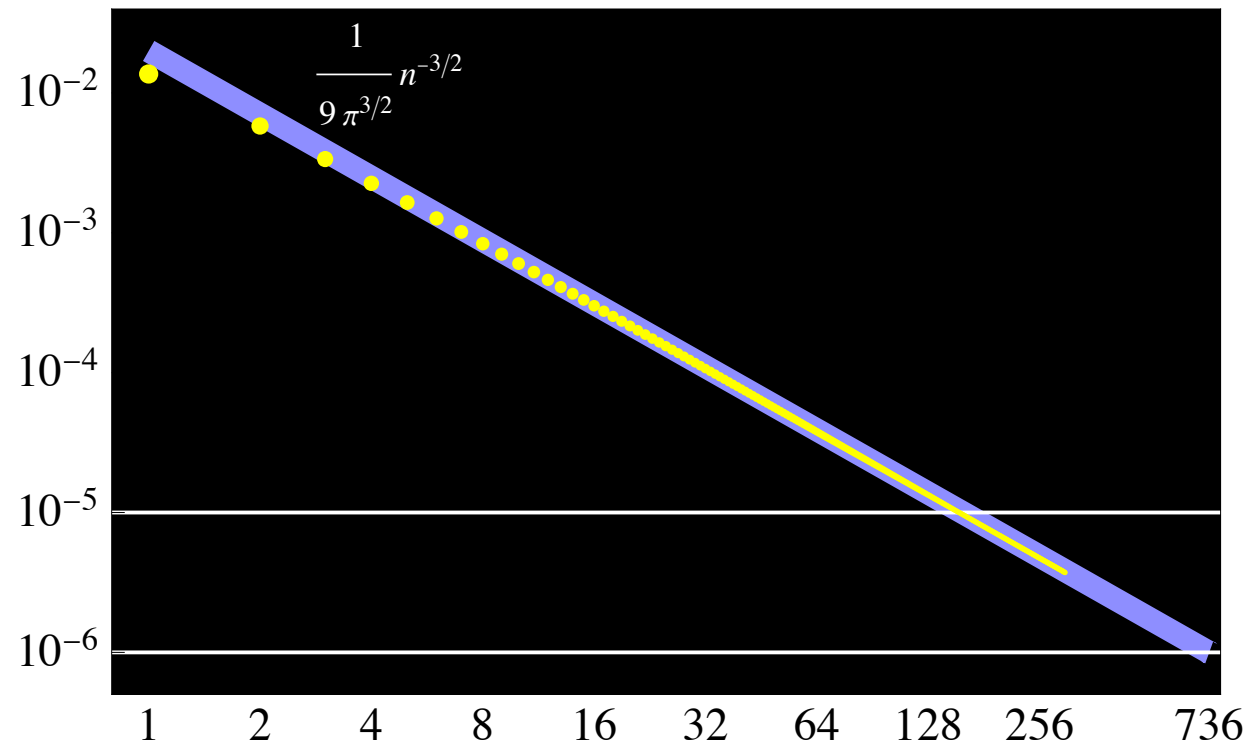

**Figure 18.** The dots mark $E_n(2\pi)$ up to $n = 300$. Computation indicates that each $E(F_n)\,(t)$ decreases as $t$ rises, so these values represent the maximum error of $F_n$. The straight line is the log–log plot of $\frac{1}{9\,\pi^{3/2}}\, n^{-3/2}$, which models the asymptotic behavior of $E_n(2\pi)$.

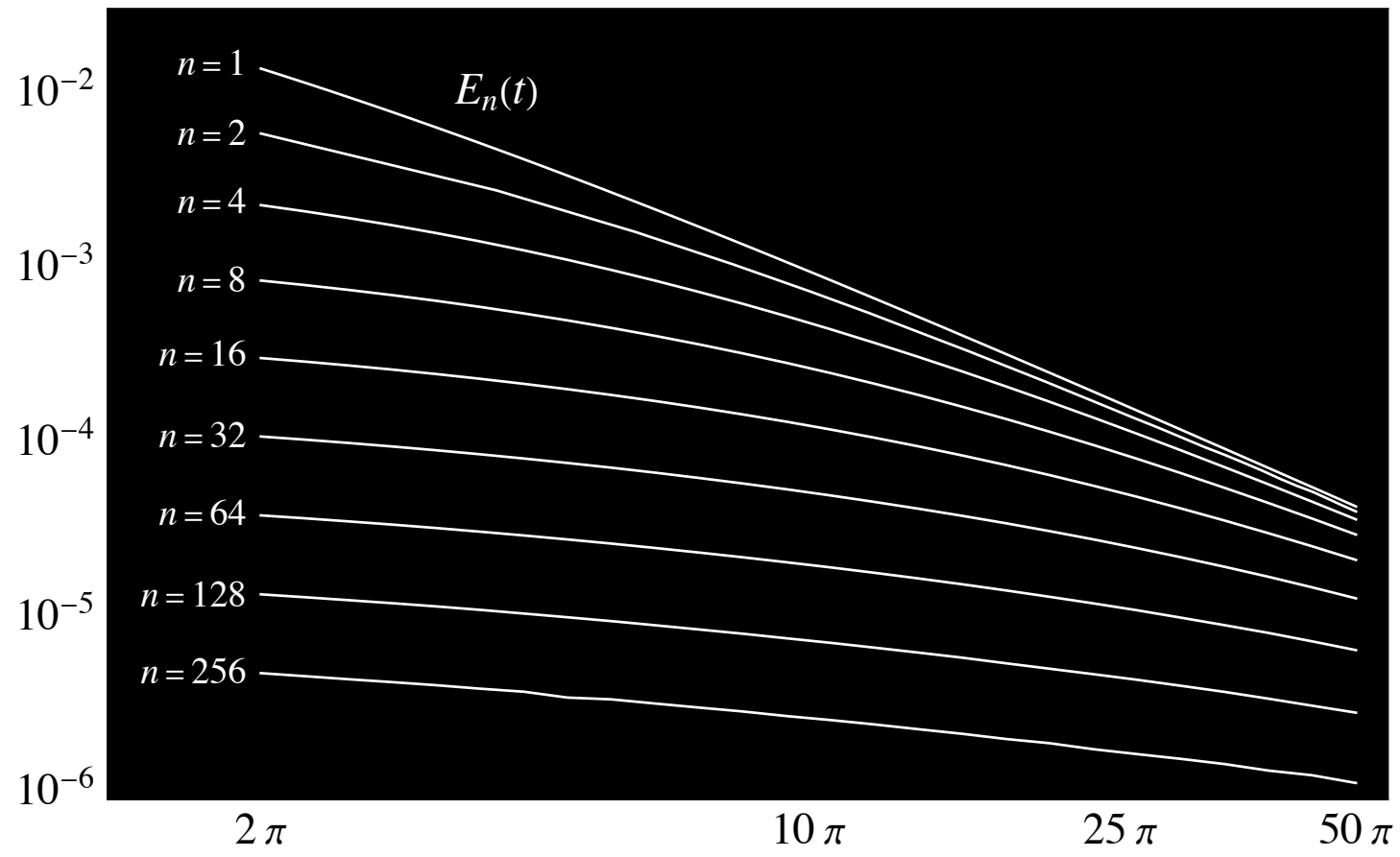

**Figure 19.** The unibike error for $F_n$ over $[2\pi, 50\pi]$. The errors all appear to decrease as $t$ increases.

Figure 20 presents evidence of convergence by checking that for each $t$, $R_{n-1}(t) - R_n(t)$ forms a Cauchy sequence. The right side of the figure compares the Cauchy differences to $E_n(t - 2\pi)$, and they appear almost identical. If their difference approaches 0, that would mean that if one goes to 0, so does the other; this indicates that the miracle of convergence is equivalent to the miracle of vanishing unibike error. Figure 21 indicates that for each $t$, the difference $(R_{n-1}(t) - R_n(t)) - E_n(t - 2\pi)$ does indeed approach 0. There is a simple reason for this relationship. Figure 22 shows that $\mu$, the closest point to $\Phi(F_n)(t)$, is close to a point on $F_n$ with the same polar angle—call it $\zeta(t)$—as $\Phi(F_n)(t)$. While we have no way to quantify the angle $\zeta(t)$, considering large values of $t$ and the spiral nature of the curve indicates that $\lim_{t\to\infty}(\zeta(t) - (t + 2\pi)) = 0$. When $n = 1$, $\zeta(t) = \alpha(t)$ and the limit is a simple consequence of $\alpha(t) = t + 2\pi + O\left(t^{-1/2}\right)$. All this provides a heuristic explanation for

(1) $\quad \lim_{n\to\infty} R_{n-1}(t + 2\pi) - R_n(t + 2\pi) - E_n(t) = 0,$

the same as the observed relationship. Example: $R_{127}(40\pi) - R_{128}(40\pi)$ is larger than $E_{128}(38\pi)$ by $1.5 \cdot 10^{-8}$. This relation plays an important role in §6 as it leads to a simple formula for a track that has much smaller unibike error than the polar square root.

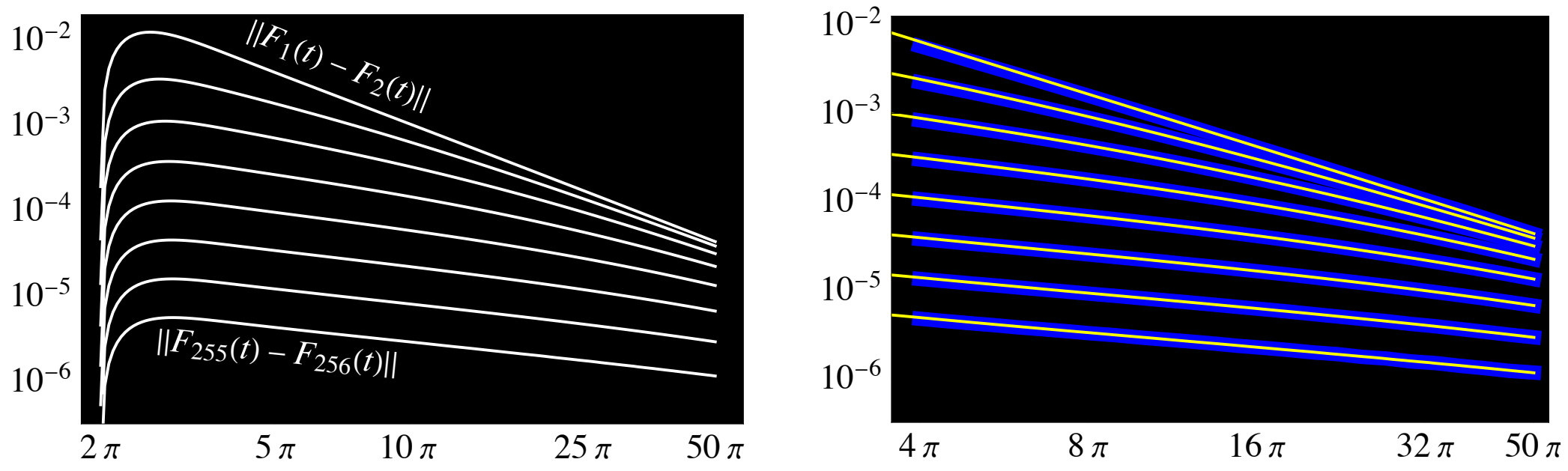

**Figure 20.** Left: $R_{n-1}(t) - R_n(t)$ for $n = 2, 4, 8, 16, 32, 64, 128, 256$. These apparently converge to 0 following a distinct pattern. The plots indicate that $(F_n)$ is uniformly Cauchy, which implies uniform convergence to a limit function. Right: The left image, with the errors $E_n(t - 2\pi)$ shown in blue (so $t$ starts at $4\pi$). This illustrates the relation (1).

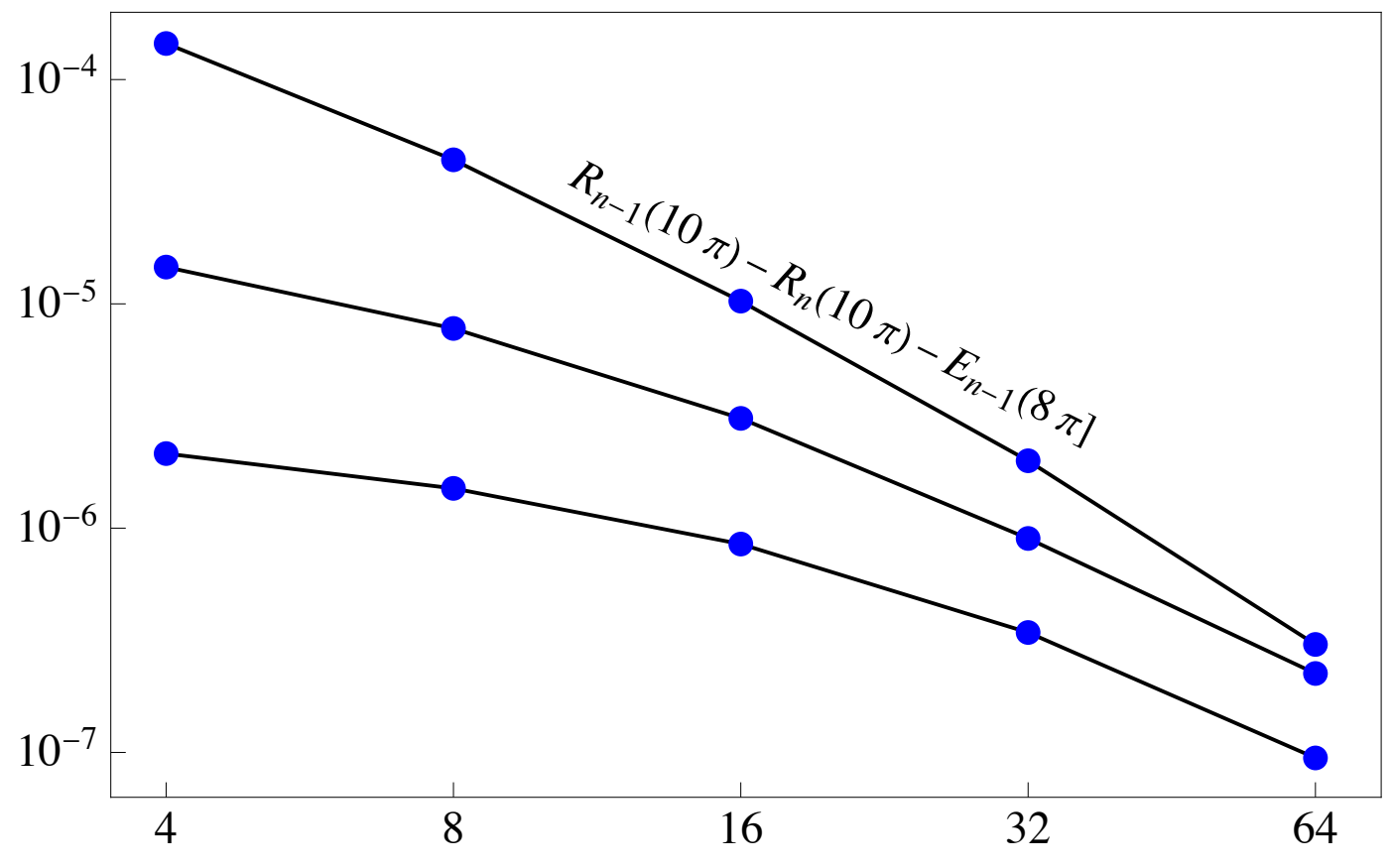

**Figure 21.** The difference $R_{n-1}(t) - R_n(t) - E_n(t - 2\pi)$ approaches 0 as $n$ rises. Here $t$ is $10\pi$, $25\pi$, and $50\pi$.

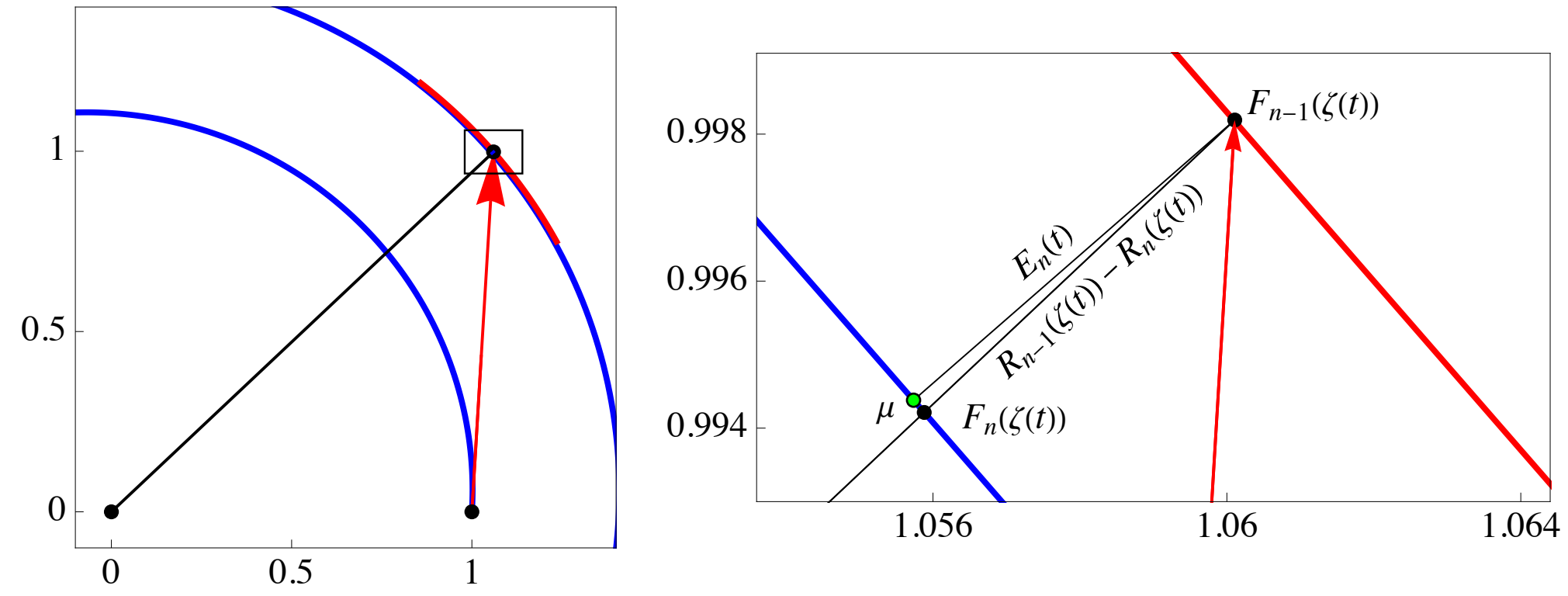

**Figure 22.** Here $n = 2$ and $t = 2\pi$ but the story applies to any $n$ and $t$. Left: $F_n$ in blue and $F_{n-1}$ in red. The front wheel point is on $F_{n-1}$. Right: Let $\zeta(t)$ be such that $\Phi(F_n)(t) = F_{n-1}(\zeta(t))$; it appears that $\zeta(t) - (t + 2\pi)$ approaches 0. The point $\mu$ is the closest point on $F_n$ to $\Phi(F_n)(t)$ and so $E_n(t)$ is the length of the line segment containing $\mu$. And $R_n(\zeta(t)) - R_{n-1}(\zeta(t))$ is the polar segment connecting the black points.

Figure 23 shows $F_{256}$, an approximation to the limit $F_\infty$ that hardly differs from the limit and so shows what a spiral unibike track looks like. The radial plot at right contains all the information for the track because it allows one to get the seed curve, and then front track iteration gives the entire curve. Figure 24 compares $\|F_1(t)\|$ with $\|F_{256}(t)\|$; there is a visible difference for small $t$, but after a small initial segment, the two functions coincide at human scale.

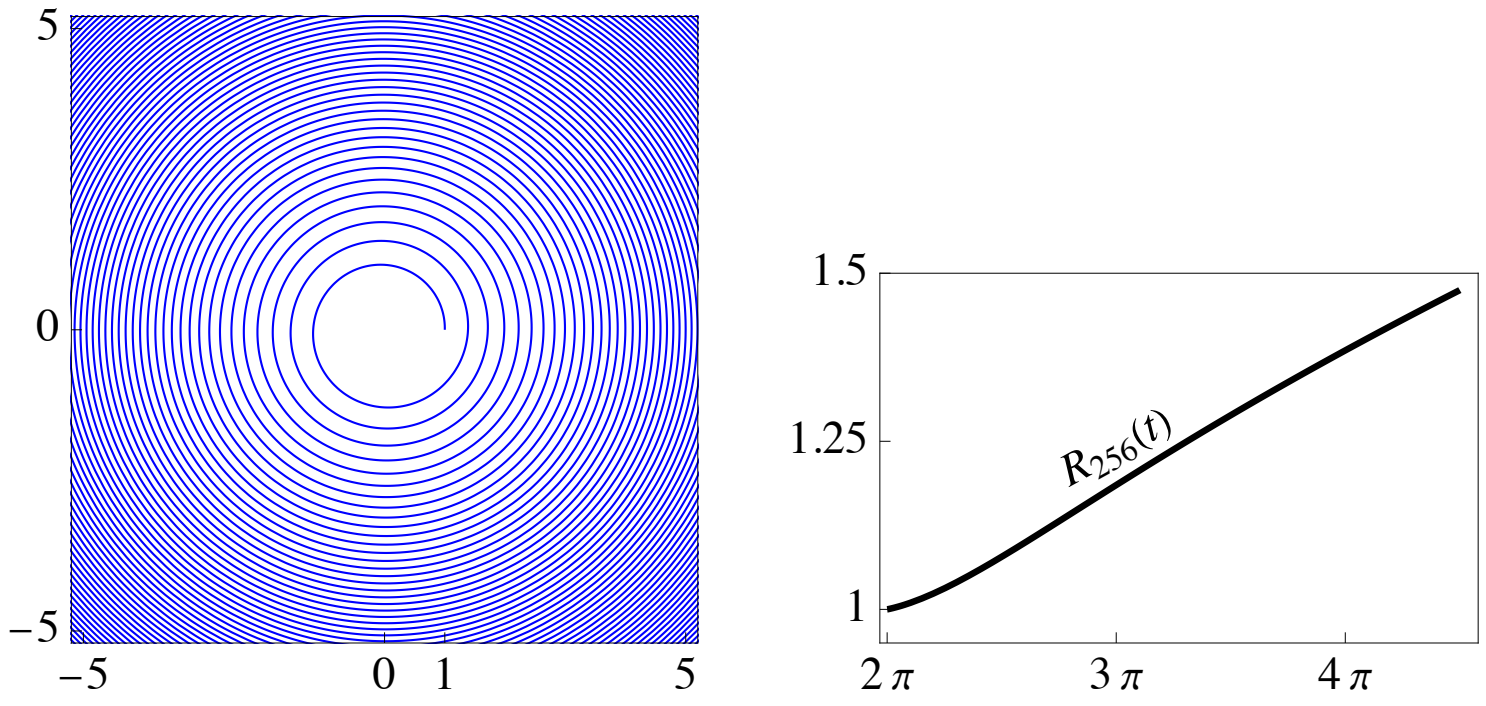

**Figure 23.** Left: $F_{256}$, which has extremely small unibike error. Right: The plot of $\|F_{256}(t)\|$ up to $4.5\pi$. The small domain is because, for the limiting curve that domain is enough to generate the full half-infinite spiral, which is presumably a unibike: just repeatedly form the front track.

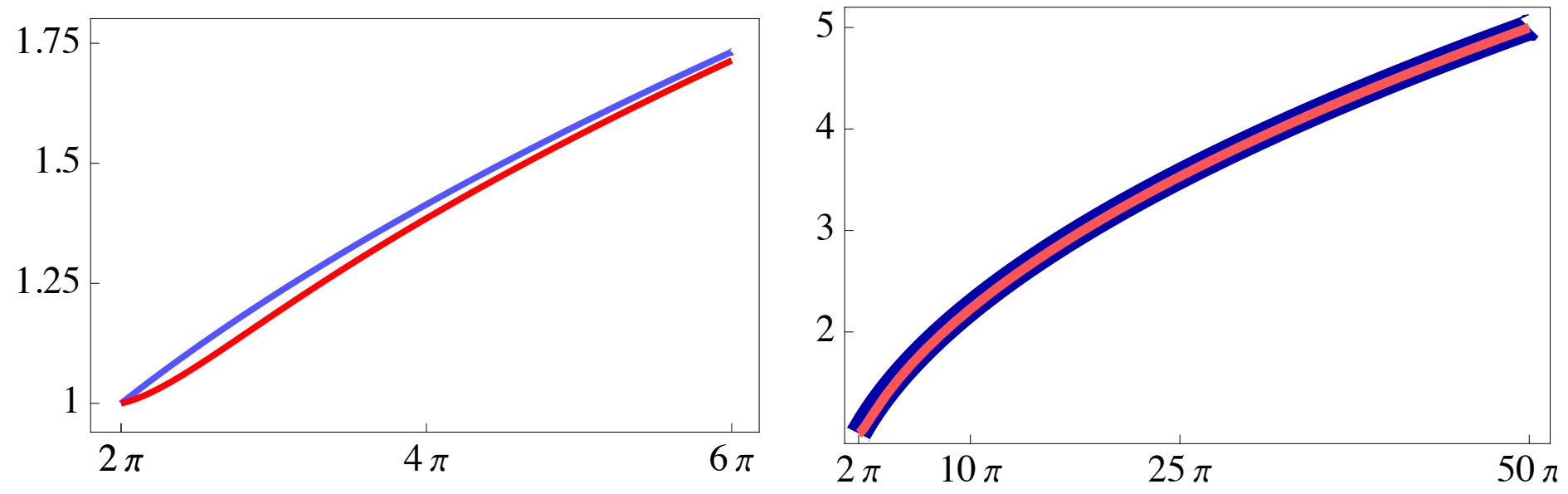

**Figure 24.** The plots compare $R_{256}(t)$ (red) with $R_1(t)$, which is $\sqrt{t/(2\pi)}$; there is some difference at the start, but overall the functions barely differ.

We can summarize these experiments in a conjecture. Perhaps this can be proved by a careful analysis of the differential equation that gives the rear track.

**Conjecture 1.** Let $(F_n)_{n\geq 1}$ be the sequence defined as follows. Start with $F_1(t) = \sqrt{t/(2\pi)}$ and define $F_n$ to be the rear track for the front track $F_{n-1}$ with the initial value $F_n(\tau) = (1, 0)$ where $F_{n-1}(\tau)$ has unit distance from $(1, 0)$ and is on the second loop of $F_{n-1}$. Reparametrize so that $t$ is the polar angle of $F_n(t)$.

(a) The sequence of functions $F_n(t)$ converges uniformly as $n \to \infty$ to an infinitely differentiable function $F_\infty(t)$ defined for $t \geq 2\pi$.

(b) $F_\infty$ is a unibike curve: the unibike error at every point $F_\infty(t)$ is 0.

(c) The curve defined by $F_\infty$ is a spiral: $F_\infty(t)$ has polar angle $t$ and $\|F_\infty(t)\|$ is monotonically increasing as $t$ rises.

One wonders whether the spiral nature of the curves would allow the rear track differential equation to be set up in just one equation, for $r'(t)$. We do not see how to do that. Also one wonders whether extrapolation can work from the sequence of radial values for fixed $t$ to predict the limiting radial value. We have had no success in that direction.

## 6. The Pochhammer Connection

Examining how the error changes as $n$ rises leads to the surprising relationship

$$(2) \qquad E_n(t) \approx E_{n-1}(t + 2\pi)\left(1 + \frac{\pi}{t+2\pi}\right),$$

relating the error at $t$ to the previous error at $t + 2\pi$. No reason for the multiplicative factor is evident, but the fit is very good, including the $\pi$ in the numerator. Figure 25 shows some data; at right, the data has been replaced by a model fit to eliminate the noise due to interpolation.

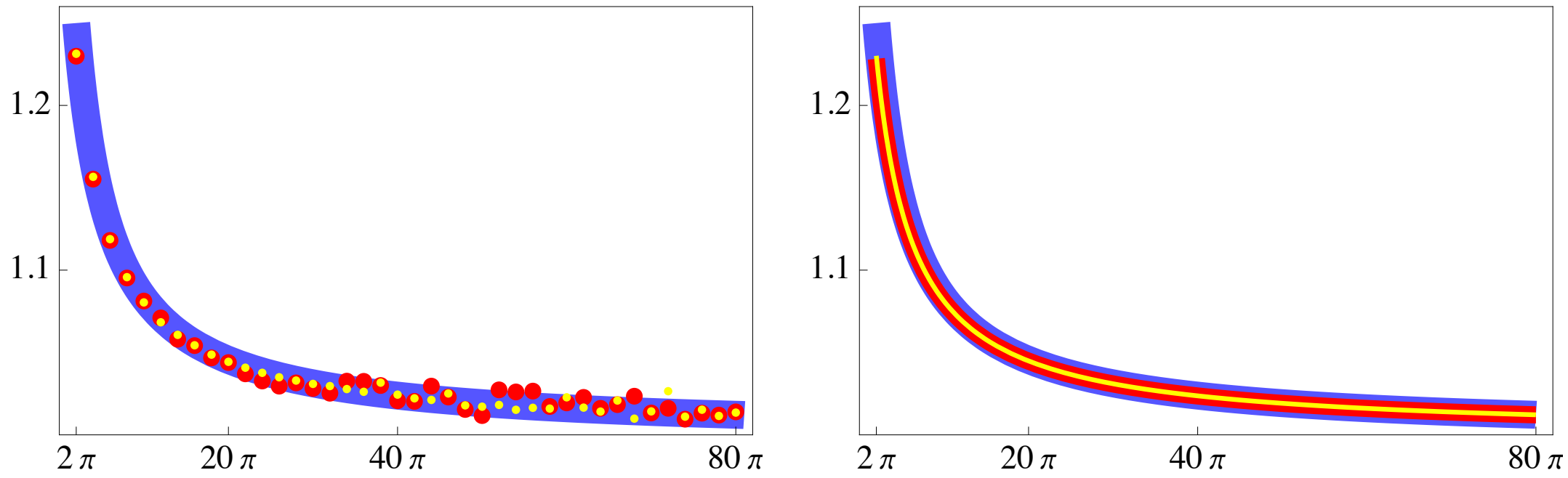

**Figure 25.** The ratio $\frac{E_n(t)}{E_{n-1}(t+2\pi)}$ compared to $1 + \frac{\pi}{t+2\pi}$ (blue), for $n = 2$ (red) and $n = 20$ (yellow). The smooth plots at right replace the noisy data for the ratios with least-squares fits to the model $1 + \frac{a}{t+b}$.

Theorem 1 proves $E_1(t) \sim (\pi/3)\, t^{-2}$, so (2) gives $E_2(t) \approx \frac{\pi}{3}(t + 2\pi)^{-2}\left(1 + \frac{\pi}{t+2\pi}\right)$ and we can repeatedly use (2) to get a finite product $E_n^*(t)$ that approximates $E_n(t)$. This product can be concisely expressed as

$$(3) \qquad E_n(t) \approx E_n^*(t) = \frac{\pi}{3}\,(2\pi(n-1) + t)^{-2}\,\frac{\left(\frac{t+3\pi}{2\pi}\right)_{n-1}}{\left(\frac{t+2\pi}{2\pi}\right)_{n-1}},$$

where $(a)_n$ is the Pochhammer function $a(a+1)\cdots(a+n-1)$. An example: $E_{64}(50\pi) = 1.3068 \cdot 10^{-5}$ and $E_{64}^*(50\pi) = 1.3036 \cdot 10^{-5}$. A proof of $E_n(t) \sim E_n^*(t)$ would imply that the unibike error goes to 0 uniformly, because $E_n^*(t)$ decreases with $t$ (a consequence of the product definition of $(a)_n$) and

$$E_n^*(2\pi) = \frac{1}{9\pi^{3/2} n^2 n!} \Gamma\left(\frac{3}{2} + n\right) = \frac{(2n+1)!!}{9\pi n^2 n! 2^{n+1}} \sim \frac{1}{9\pi^{3/2}} n^{-3/2}.$$

If one wants to know only that the error goes to 0 (assuming (2) is asymptotically true), then this works:

$$E_n^*(2\pi) = \frac{1}{9\pi^{3/2} n^2 n!} \Gamma\left(\frac{3}{2} + n\right) \le \frac{1}{9\pi^{3/2}} \frac{(n+1)!}{n^2 n!} = \frac{1}{9\pi^{3/2}} \frac{n+1}{n^2} \le \frac{1}{n}.$$

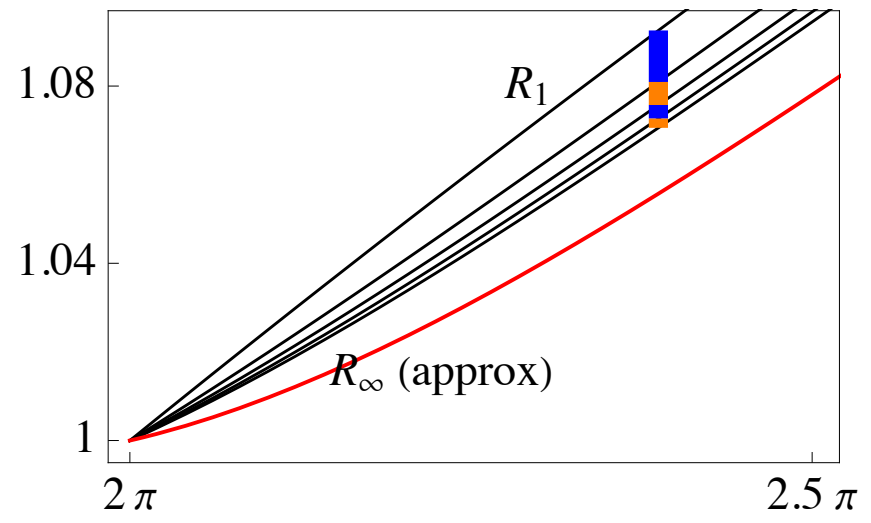

**Figure 26.** The differences $R_n - R_{n+1}$, when subtracted from $R_1$, give $R_\infty$.

For each $t$, the polar radius $R_n(t)$ appears to decrease to a limit $R_\infty(t)$ (see Fig. 26), so we can use (1) and (3) to approximate $R_\infty$ with a hypergeometric function.

$$R_\infty(t) = R_1(t) - (R_1(t) - R_2(t)) - (R_2(t) - R_3(t)) - \ldots \approx \sqrt{\frac{t}{2\pi}} - (E_2(t-2\pi) + E_3(t-2\pi) + \ldots) \text{ (by (1))}$$
$$\approx \sqrt{\frac{t}{2\pi}} - (E_2^*(t-2\pi) + E_3^*(t-2\pi) + \ldots) \text{ (by (3))}$$
$$= \sqrt{\frac{t}{2\pi}} - \frac{\pi}{3t^3}(t+\pi)\, {}_4F_3\left(1, \frac{t}{2\pi}, \frac{t}{2\pi}, \frac{t}{2\pi} + \frac{3}{2}; 1 + \frac{t}{2\pi}, 1 + \frac{t}{2\pi}, 1 + \frac{t}{2\pi}; 1\right),$$

where ${}_4F_3$ is a generalized hypergeometric function ${}_pF_q$. Because ${}_4F_3$ is defined as a sum of Pochhammer values, we are not really getting anything new, but this gives us a nice unibike approximation using no differential equations. Define $H_1(t) = \left(\sqrt{\frac{t}{2\pi}} - \frac{\pi}{3t^3}(t+\pi)\, {}_4F_3\left(1, \frac{t}{2\pi}, \frac{t}{2\pi}, \frac{t}{2\pi} + \frac{3}{2}; 1 + \frac{t}{2\pi}, 1 + \frac{t}{2\pi}, 1 + \frac{t}{2\pi}; 1\right)\right)(\cos t, \sin t)$. It appears that the subtracted term, $\sqrt{\frac{t}{2\pi}} - \|H_1(t)\|$, is $\frac{1}{3}t^{-1} + \frac{\pi}{3}t^{-2} + O(t^{-3})$, which means that $H_1(t) \sim F_1(t)$, but for small $t$ they differ by a lot (e.g., $H_1(2\pi) = (0.91\ldots, 0)$ compared to $(1, 0)$ for $F_1$).

It is possible, but tedious, to compute rear tracks starting with $H_1$, as was done in §5 starting from $F_1$. The difficulty is the evaluation of both $H_1(t)$ and $H_1'(t)$. One can spend several hours getting an interpolating function based on a fine mesh of values of $H_1$ and then using the interpolant to generate the rear track sequence. But a niftier approach is to use the observation above in the form $\|H_1(t)\| \approx \sqrt{\frac{t}{2\pi}} - \frac{1}{3t}$. There has been much research into the asymptotics of hypergeometric functions, but this result appears to have not been proved. Numerical computations give solid support, so we can conjecture as follows, where the form is simplified by replacing $t$ by $2\pi t$ and eliminating irrelevant coefficients.

**Conjecture 2.** ${}_4F_3\left(1, t, t, t + \frac{3}{2}; t+1, t+1, t+1; 1\right) \sim 2t^2$.

This conjectured asymptote to $H_1$ leads to an adjustment of the polar square root to $G_1(t) = \left(\sqrt{\frac{t}{2\pi}} - \frac{1}{3}t^{-1}\right)(\cos t, \sin t)$, a nicely simple formula that is a significant improvement over $F_1$ as a unibike. And $G_1$ does define a spiral: the polar radius increases with $t$. In Theorem 2 we will prove that the unibike error for $G_1$ is much smaller than $(\pi/3)t^{-2}$, the asymptotic error for $F_1$. Figure 27 compares the unibike errors for $F_1$ and $G_1$, and the correct model for $E(G_1)$ is evidently about $0.63\,t^{-5/2}$; the coefficient appears to be $\pi/5$, but it might not be exactly this. The graph also shows the error for $K_1$, which uses one more term: $\|K_1\| = \sqrt{\frac{t}{2\pi} - \frac{1}{3}t^{-1} - \frac{\pi}{3}t^{-2}}$. The maximum unibike errors for $H_1$, $G_1$, and $K_1$ are $1/531$, $1/126$, and $1/523$, respectively, so $K_1$ is a fine proxy for $H_1$. The error for $H_1$ is also shown in Figure 27; it is very different from $E(K_1)$.

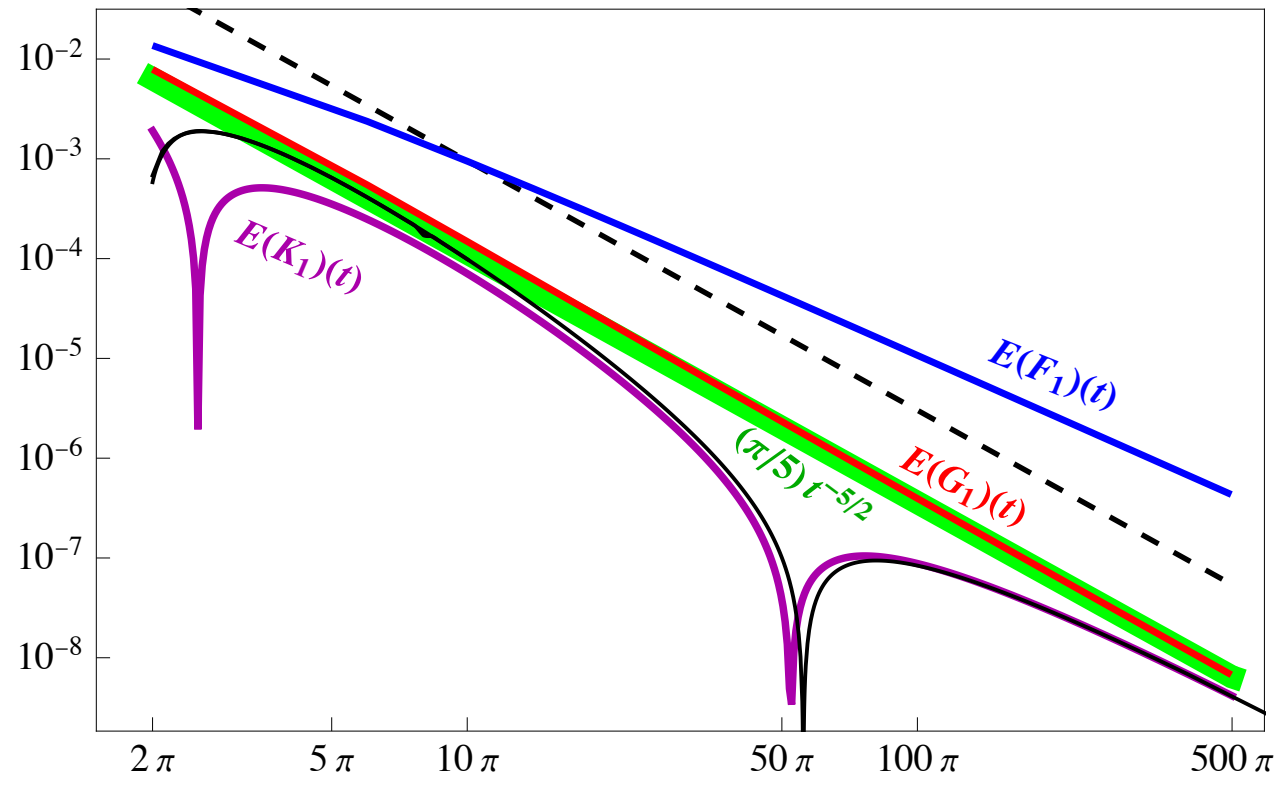

**Figure 27.** Log–log plots of the unibike error for $F_1$ (blue) and $G_1$ (red), with a fitting curve for $E(G_1)$ shown in green; the coefficient appears to be close to $\pi/5$. The dashed line is the bound on the proof of Theorem 2, with coefficient about 5.28, compared to 0.628 in the better fit (green). The purple graph shows the error for $K_1$, which slightly improves on $G_1$ (the maximum errors for these two curves are $1/523$ vs. $1/126$). The two dips to 0 indicate places where the front curve for $K_1$ intersects $K_1$. The black graph is the error for $H_1$; a front–rear crossing occurs near $t = 56\pi$.

**Theorem 2.** The error $E(G_1)(t)$ is asymptotically bounded by $6\,t^{-5/2}$.

**Proof.** Let $G_0 = \Phi(G_1)$, the front track. Any point $G_1(s)$ provides an upper bound on the error: the minimum distance from $G_0(t)$ to $G_1(s)$. Let $s = \alpha(t) + \frac{2}{3}\pi\,t^{-2}$, where $\alpha(t)$ is as in §2. The quadratic term in $s$ was found by numerical work using the location of the true closest point, and identifying the constant $2\pi/3$ in a least-squares fit to the $s - \alpha(t)$. So now it remains only to find the asymptotic value of $\|G_0(t) - G_1(s)\|$. This is a complicated expression but we can get the desired result with some computer assistance when we use the squared distance $\delta = (G_0(t) - G_1(s)) \cdot (G_0(t) - G_1(s))$. Then *Mathematica* gets the result via $\lim_{t\to\infty} \delta\, t^5 = (\pi/72)\left(9 + 64\,\pi^2\right) < 6^2$. Code for this last step is in the Appendix. □

Knowing that $G_1$ is an improvement over $F_1$, we can hope for nicely reduced error when we compute the sequence of rear tracks from $G_1$ as in §5. The results are very good. Each $E(G_n)(t)$ appears to be monotonically decreasing. And as $n$ increases the error decreases similarly to $E(F_n)$, but with a much smaller coefficient (Fig. 28). For $E(F_n)$ the maximum error is about $0.02\,n^{-3/2}$ while $E(G_n)$ is asymptotic to about $0.00033\,n^{-3/2}$. The unibike error for $G_{72}$ is under $10^{-6}$, a bound not attainable by computation for the polar square root, where $F_{736}$ would be required. Figure 28 includes the error $E(K_n)(2\pi)$, where $K_n$ is the sequence of rear tracks starting from $K_1$. For some small values of $n$ the error at $2\pi$ is not the maximum error (e.g., $K_{26}$), but eventually the normal pattern of decreasing error returns. So $K_{300}$ is quite good, with maximum error under $7 \cdot 10^{-8}$. The steep drop in Figure 28 at $n = 26$ is because of a crossing of the front and rear tracks that occurs near $2\pi$, which leads to near-zero error at $2\pi$. Such crossings exist already for $K_1$ (two are evident in Fig. 27), but they move left as $n$ rises and eventually disappear. Figure 29 shows the data, which is very noisy but does support the hypothesis the error for the true $K_{128}$ path only decreases as $t$ rises.

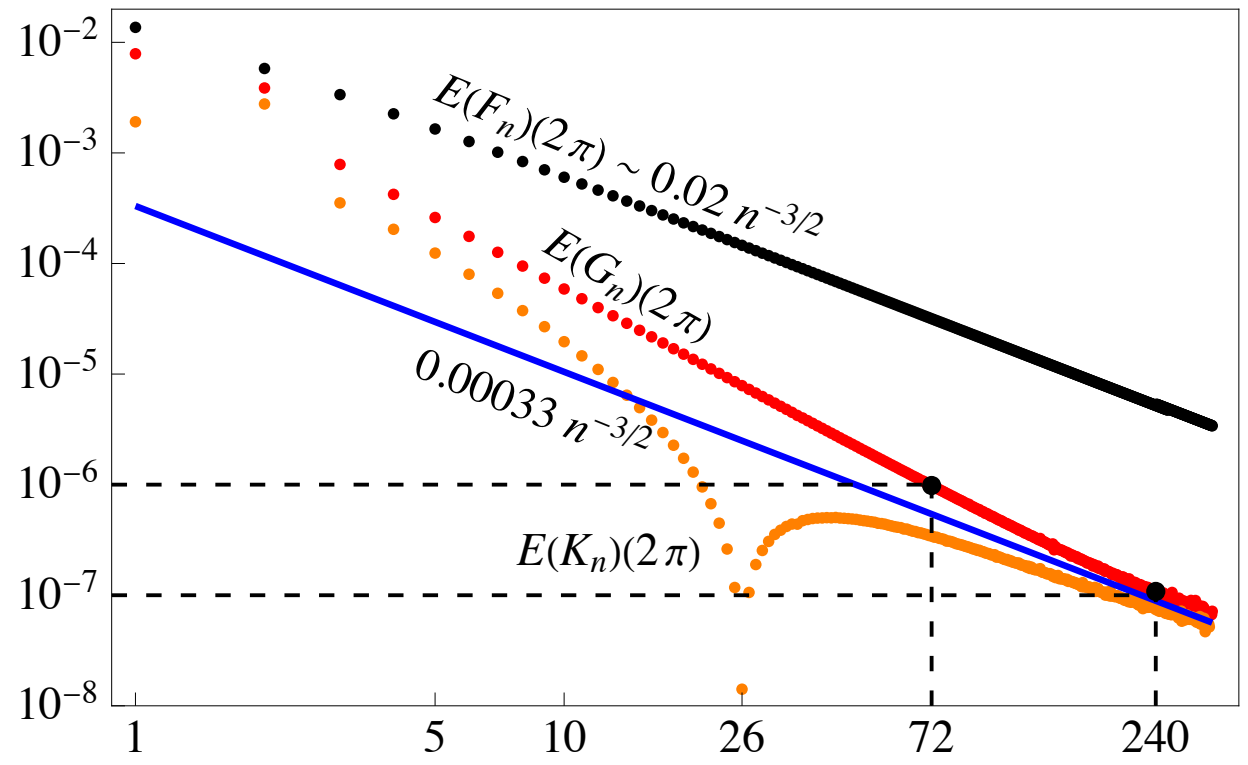

**Figure 28.** The maximum unibike error for $G_n$ (red) is apparently asymptotic to about $\frac{1}{3000}\, n^{-3/2}$; this is similar to the behavior of $E(F_n)$, but with a much smaller constant. The maximum error of $K_n$ (the orange graph) is quite a bit smaller than $\max(E(G_n))$, but the asymptotic behavior seems unchanged.

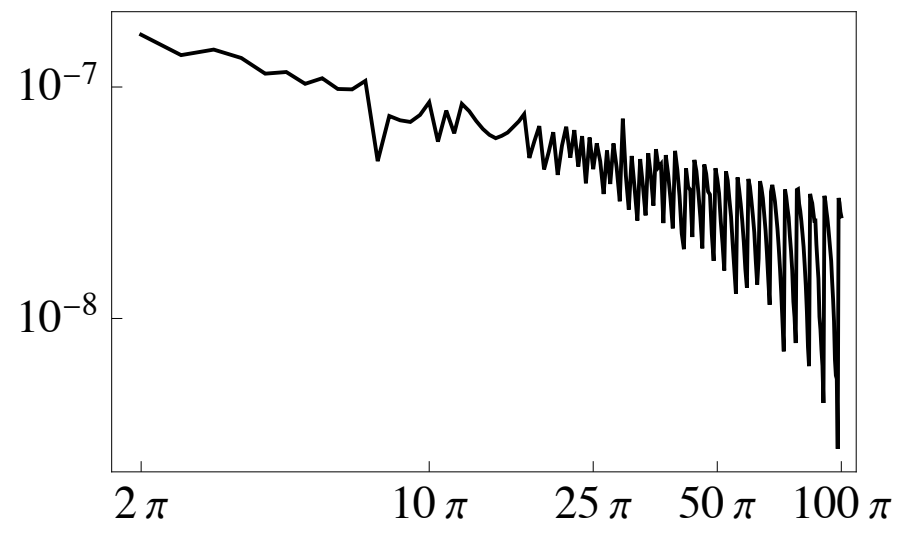

**Figure 29.** The error for $K_{128}$ is very noisy, but it appears that the ideal version of this function is decreasing from its initial value at $t = 2\,\pi$.

The various improvements here are quite large. The maximum polar square root unibike error is about 1 part in 73; the same for $K_{320}$ is less than 1 part in 19,000,000.

A final idea is to consider a linear combinations of $S_n$ and $S_{2\,n}$, where $S$ is any of $F$, $G$, or $K$. This is inspired by the classic idea used by Huyghens to improve Archimedes's approximation to $\pi$; more generally it is known as Romberg's method or Richardson extrapolation. Our setup is not obviously one where the method applies, but it is easy enough to try various weighted averages, and indeed, we get significant improvement. With $\sigma = 6.83$, consider $K_{12}^{*} = (\sigma\, K_{12} - K_6) / (\sigma - 1)$. Then $E(K_{12}^{*})$ is bounded by $10^{-6}$ (Fig. 30) so this breaks the one-part-in-a-million goal using only 11 differential equation solutions. This is the spiral shown in Figure 4.

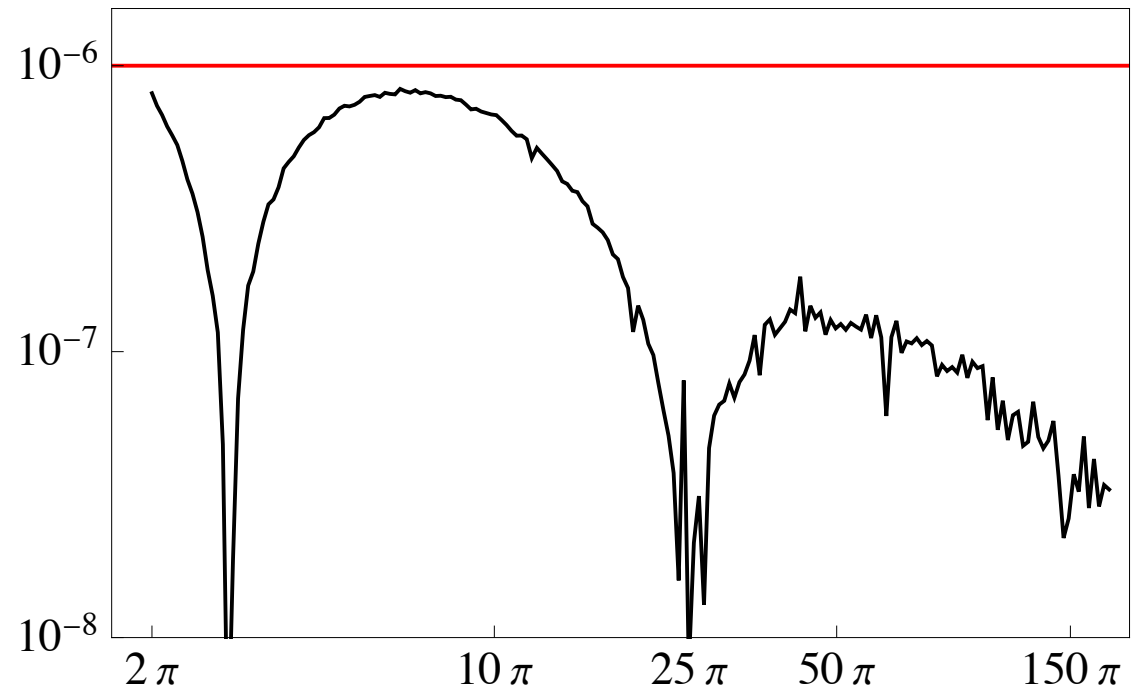

**Figure 30.** The unibike error for $K_{12}^{*}$ is bounded by $10^{-6}$, thus achieving one-part-in-a-million error with only 11 differential equations.

## 7. Conclusion

Finn's example is very complicated and it is natural to think that there can be a unibike curve in the shape of a simple spiral. The very simple polar square root curve has maximum unibike error of 1 part in 73. Then the rear track iteration leads to a curve with error that is at most 1 part in 268,000. An adjustment to the polar square root gets to 1 part in 14 million, and a final enhancement leads to the spiral curve $K_{320}$ with unibike error less than one part in 19 million. There is a pattern to the error decrease, so all these methods provide evidence that there is a limiting spiral that is a perfect unibike. Perhaps one can carefully analyze the differential equation to discover a rigorous explanation of why moving to the rear track decreases the error as the computations indicate, and so get a rigorous proof that the limiting curve exists and is a perfect spiral unibike.

The author is grateful to Bob Bixler, Michael Elgersma, Antonín Slavík, and Michael Trott for helpful comments on the unibike problem.

# Appendix of *Mathematica* Code

## Code for asymptotics in proof of Theorem 1(a)

```
F[t_] := Sqrt[t / (2 π)] {Cos[t], Sin[t]}; Fder = Normalize[D[F[t], t]];
FFront[t_] := F[t] + Fder; α[t_] := 2 π + t + ArcTan[2 t / (1 + (Sqrt[t] Sqrt[1+4 t^2]) / Sqrt[2 π])];
Series[Simplify[Norm[FFront[t]] - Norm[F[α[t]]], t ≥ 2 π], {t, ∞, 4}]
```

$$\frac{\pi}{3\,t^2}+\frac{1}{2}\sqrt{\frac{\pi}{2}}\left(\frac{1}{t}\right)^{5/2}-\frac{11\,\pi^2}{15\,t^3}-\frac{11\,\pi^{3/2}\left(\frac{1}{t}\right)^{7/2}}{6\,\sqrt{2}}+\frac{\pi\left(-175+412\,\pi^2\right)}{280\,t^4}+O\left[\frac{1}{t}\right]^{9/2}$$

## Code for proof of Theorem 1(b).

```
left[{A_, B_}, C_] := (A - C).({1, -1} Reverse[B - C]); Q = F[α[t]];
{Limit[TrigExpand[t left[{Q, Q - F′[α[t]]}, F[t + 4 π]]], t → ∞],
 N[t left[{Q, Q - F′[α[t]]}, F[t + 4 π]] /. t → 10^20 π, 20]}
```

{-π, -3.1415926534955123342}

## Complete code to generate the rear track for the polar square root

```
F[t_] := Sqrt[t / (2 π)] {Cos[t], Sin[t]}; G[t_] := {Gx[t], Gy[t]}; b = F[t] - G[t];
t0 = SolveValues[EuclideanDistance[F[t], {1, 0}] == 1 && 4 π < t < 5 π, t][[1]];
GDE = NDSolveValue[{G′[t] == (F′[t].b) b, G[t0] == {1, 0}}, G[t], {t, t0, 60 π}];
ParametricPlot[{F[t], GDE}, {t, t0, 10 π},
 PlotStyle → {{Thickness[0.025], Blue}, {Medium, Red}}]
```

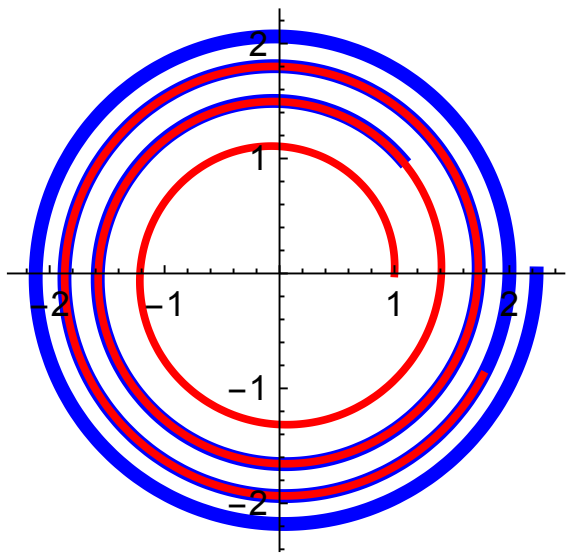

## Code for the limit in the proof of Theorem 2

```
α[t_] := 2 π + t + ArcTan[(4 Sqrt[π] t) / (2 Sqrt[π] + Sqrt[2] Sqrt[t] Sqrt[1 + 4 t^2])];
G[t_] := (Sqrt[t] / Sqrt[2 π] - 1 / (3 t)) {Cos[t], Sin[t]};
Gd = Simplify[Normalize[∂_t G[t]], t > 2 π];
δ = TrigExpand[Simplify[G[t] + Gd - G[α[t] + (2 / 3) π t^-2]]];
Together[Sqrt[Limit[δ.δ t^5, t → ∞]]]
```

$$\frac{1}{6}\sqrt{\frac{1}{2}\,\pi\left(9+64\,\pi^2\right)}$$